\documentclass[12pt]{article}
\usepackage{color}
\usepackage{amscd,amsfonts,amssymb,amscd,amsmath,latexsym}
\usepackage[hypertex,backref]{hyperref}
\usepackage[all]{xy}
\makeatletter
\renewcommand{\subsubsection}{\@startsection
{subsubsection}
{1}
{0mm}
{0mm}
{0mm}
{\normalfont\normalsize\itshape}}
\makeatother
  
\usepackage{mathrsfs}

\title{The $f$-invariant and index theory}
\author{Ulrich Bunke\thanks{NWF I - Mathematik,
Universit{\"a}t Regensburg,
93040 Regensburg,
GERMANY, ulrich.bunke@mathematik.uni-regensburg.de} \and Niko Naumann \thanks{NWF I - Mathematik,
Universit{\"a}t Regensburg,
93040 Regensburg,
GERMANY, niko.naumann@mathematik.uni-regensburg.de} }
 
\newtheorem{theorem}{Theorem}[section] 
\newtheorem{prop}[theorem]{Proposition}
\newtheorem{lem}[theorem]{Lemma}

\newtheorem{ddd}[theorem]{Definition}

\renewcommand{\lim}{{\tt lim\:}}

\newcommand{\Z}{\mathbb{Z}}

\newcommand{\proof}{{\it Proof.$\:\:\:\:$}}

\newcommand{\R}{\mathbb{R}}

\newcommand{\Q}{\mathbb{Q}}

\renewcommand{\det}{{\tt det}}
\renewcommand{\sinh}{\mathrm{sinh}}

\newcommand{\C}{\mathbb{C}}

\newcommand{\Sf}{{\tt  Sf}}

\newcommand{\cZ}{\mathcal{Z}}

\newcommand{\cM}{\mathcal{M}}

\newcommand{\cE}{\mathcal{E}}

\newcommand{\cY}{\mathcal{Y}}

\newcommand{\Ext}{{\tt Ext}}

\newcommand{\cD}{\mathcal{D}}

\newcommand{\Imm}{{\tt Im}}

\newcommand{\ee}{{\tt e}}

\newcommand{\cN}{\mathcal{N}}

\newcommand{\id}{{\tt id}}

\newcommand{\nat}{\mathbb{N}}

\newcommand{\Td}{{\mathbf{Td}}}

\def\imath{{i}}

\def\hB{\hspace*{\fill}$\Box$ \newline\noindent}

\newcommand{\ind}{{\tt index}}

\def\hB{\hspace*{\fill}$\Box$ \\[0.5cm]\noindent}

\newcommand{\pr}{{\tt pr}}

\newcommand{\ch}{{\mathbf{ch}}}

\newcommand{\hA}{\hat{\mathbf{A}}}

\newcommand{\ZN}{{{}^N\Z}}
 \newcommand{\Dirac}{{\not\hspace{-1.3mm}\mathcal \cD}}

\begin{document}

\maketitle
 
\begin{abstract}
In this paper we prove a tertiary index theorem which relates
a spectral geometric and a homotopy theoretic invariant 
of an almost complex manifold with framed boundary.
It is derived from the index theoretic and homotopy theoretic
versions of a complex elliptic genus and interestingly related with the structure
of the stable homotopy groups of spheres.
\end{abstract}

\tableofcontents

\section{Introduction}\label{uiwrwerwrwer}

The archetypical assertion in index theory is an equality 
\begin{equation}\label{dwdwqdwqd}
\ind^{an}=\ind^{top}
\end{equation} of a topological and analytical index. To be more specific, we consider the Dirac operator $\Dirac_M$ on a closed almost complex manifold $M$ of dimension $2k$.  In order to define this operator we must choose a Riemannian metric
and a $Spin^c$-extension of the Levi-Civita connection to the $Spin^c$-principal bundle determined by the almost complex structure. The index
$\ind(\Dirac_M)\in \Z$ of $\Dirac_M$ is defined as the super dimension of its   kernel. 
It is independent of the choice of the geometric structures and actually only depends on the almost-complex bordism class $[M]\in MU_{2k}$. In this way the analytical index gives a homomorphism
$$\ind^{an}:MU_{2k}\to \Z\ ,\quad \ind^{an}([M]):=\ind(\Dirac_M)\ .$$

Complex $K$-theory is a complex oriented generalized cohomology theory. The complex  orientation is a map of spectra $$\theta:MU\to K\ .$$
On coefficients it induces the topological index homomorphism
$$\ind^{top}:MU_{2k}\to K_{2k}\cong \Z\\ ,\quad \ind^{top}([M]):=\theta_{2k}([M])  .$$
The equality (\ref{dwdwqdwqd}) is then a special case of the Atiyah-Singer index theorem \cite{MR0236950}.

The equality (\ref{dwdwqdwqd}) 
is the primary index theorem. The  main purpose of the present paper is to pursue a method 
to construct higher derived topological and analytical index quantities
and to prove their equality.  As it turns out our example is very interestingly related to the stable homotopy groups of spheres. 
The present paper gives the first example of a tertiary index theorem.

Let us explain the rough idea right now. We start with the secondary invariants. Their construction depends on the fact that both, the topological and analytical index are almost local.
More precisely,
the topological index can be calculated
as an evaluation $ \langle \Td(TM),[M]\rangle$  of a characteristic class of the almost complex  tangent bundle of $M$. 
 Assume that we cut the manifold $M$ in halfs along a hyper surface $N$,
$M=M_0\cup_N M_1$, and that the tangent bundle $TN$ is trivialized (framed)
as a (stable) almost complex bundle. Then we can refine 
 the Todd class to a rational relative cohomology class so that
\begin{equation}\label{ufwiwfwefewfwefwefwef}
\langle \Td(TM_0,N),[M_0,N]\rangle+ \langle \Td(TM_1,N),[M_1,N]\rangle= \langle \Td(TM),[M]\rangle\ .
\end{equation} 
In this fomula the integer on the right-hand side is expressed as the sum
of two rational numbers. It follows that the class
$[\langle \Td(TM_0,N),[M_0,N]\rangle]_{\R/\Z}\in \R/\Z$ only depends on the  framed bordism class
$[N]\in S_{2k-1}$. In this way we get the secondary topological index 
\begin{equation}\label{uifuewfwefewf2}e^{top}:S_{2k-1}\to \R/\Z\ , \quad e^{top}([N]):=[\langle \Td(TM_0,N),[M_0,N]\rangle]_{\R/\Z}\ .\end{equation}
Here $k\ge 1$, $S$ denotes the sphere spectrum,  and $S_{2k-1}=\pi^s_{2k-1}(S^0)\cong \Omega^{fr}_{2k-1}$ is the $2k-1$'th  stable homotopy group of
$S^0$ which can be identified with the corresponding framed bordism group by the Pontrjagin-Thom construction. 
The notation $e^{top}$ is not accidential since this it in fact the famous $e$-invariant introduced in \cite{adams}. 

The almost locality of the analytical index can be expressed in the fact, that one can formulate suitable boundary conditions in  order to define Fredholm operators
$\Dirac_{M_i}$ whose indices sum up to $\ind(\Dirac_M)$. The choice of the boundary condition on the analytic side is a refinement of the relative $K$-homology class $[\Dirac_{M_0}]\in K_{-2k}(M_0,N)$ to an absolute class in $K_{-2k}(M_0)$.
It corresponds to the refinement of
$\Td(TM_0)\cap [M_0]\in H_{*}(M_0,N;\Q)$ to the class
$\Td(TM_0,N)\cap [M_0,N]\in 
H_*(M_0;Q)$. In the present paper we consider boundary conditions of Atiyah-Patodi-Singer type. In fact, the analysis of the boundary contribution to the index formula led  \cite[Theorem 4.14]{MR0397798} to define the 
analytic secondary index
\begin{equation}\label{uifuewfwefewf}
e^{an}:S_{2k-1}\to \R/\Z\ .
\end{equation}
The details will be explained in Section \ref{duwiqdqwd}, in particular see (\ref{zuzu2}).

The secondary index theorem states
$$e^{an}= e^{top}\ .$$
An obvious advantage of the analytic formula (\ref{zuzu2}) for 
$e^{an}([N])$ is that in contrast to the topological expression (\ref{uifuewfwefewf2}) it is intrinsic in $N$. This fact has fruitfully been exploited in  \cite{deninger} as will be explained in greater detail below.

The idea of the construction of tertiary invariants is essentially to apply the constructions above to $e^{an}$ and $e^{top}$ in place of $\ind^{an}$ and $\ind^{top}$, respectively.  This is not a canonical matter but involves choices, e.g.  as a first step one must extend the definition of the $e$-invariant to almost complex manifolds instead of framed ones.  In the present paper we  choose to work with Dirac operators related with complex elliptic genera. Another example using the action of Adams operations will be discussed in a \textcolor{black}{subsequent} paper \textcolor{black}{(in preperation).}

Roughly, the Dirac operator associated to the complex elliptic genus is the twisted operator $\Dirac_M\otimes C(TM)(q)$, where
$C(TM)(q)$ is a certain formal power series of bundles (\ref{zuidwqdqwdwqdqwd}) derived from
the almost complex tangent bundle.  For the purpose of this introduction
just note that $\ind(\Dirac_M\otimes C(TM)(q))\in E^\Gamma_{2k}[[q]]\subset {}^N\Z[[q]]$
is a formal power series which is the $q$-expansion of an integral modular form. The exact notation will be explained in Section \ref{sdvdsv}.
In this case the primary invariant is a homomorphism 
$$\ind:MU_{2k}\to \tilde E^\Gamma_{2k}$$ having values in the coefficients of a complex-oriented
elliptic cohomology theory $\tilde E^\Gamma$.  We now consider a partition $M=M_0\cup_NM_1$ along a not necessarily framed manifold $N$.
The boundary contribution to the index
theorem for an appropriate Fredholm extension of $\Dirac_{M_0}\otimes C(TM_0)(q)$ is the $\eta$-invariant
of $\Dirac_N\otimes C(TN)(q)$. Since it represents the analog of the $e$-invariant above we denote it by
$e_{ell}(N)$ for the moment.
 It follows from the APS-index theorem  \cite{MR0397798} that 
$$e_{ell}(N)\in E^\Gamma_{\C,2k}[[q]]+ {}^N\Z [[q]]+\C\subset \C[[q]]\ .$$
This fact can be considered as the analog of the integrality of the ordinary index.

  In order to construct tertiary invariants we now proceed as above.
We consider a partition $N:=N_0\cup_ZN_1$ along a hyper surface $Z$ whose (stably) almost complex tangent bundle is trivialized. Instead of the index of a boundary value problem we consider the $\eta$-invariant of an appropriate boundary value problem which we denote by 
$e_{ell}(N_0,Z)\in  \C[[q]]$ for the moment. 
Almost locality of $e_{ell}$ manifests itself in
the equality
\begin{equation}\label{zdqdqduqdzwqdqwdqwdqwdqwd}
e_{ell}(N_0,Z)+e_{ell}(N_1,Z)=e_{ell}(N)
\end{equation}
which is the analog of (\ref{ufwiwfwefewfwefwefwef}).
While (\ref{ufwiwfwefewfwefwefwef}) is a consequence of the APS-index theorem
for manifolds with boundary  \cite{MR0397798} the proof of (\ref{zdqdqduqdzwqdqwdqwdqwdqwd}) employs in a similar manner the more recent index theorem  \cite{math.DG/0201112} for manifolds with corners.
In (\ref{zdqdqduqdzwqdqwdqwdqwdqwd})  the element of $E^\Gamma_{\C,2k}[[q]]+ {}^N\Z [[q]]+\C$ on the right-hand side is expressed as a sum of two elements of $\C[[q]]$.
This easily implies that the class
$$[e_{ell}(N_0,Z)]\in  \frac{\C[[q]]}{E^\Gamma_{\C2,k}[[q]] + {}^N\Z [[q]]+\C}$$
only depends on the framed bordism class $[Z]\in S_{2k-2}$.
In this way we define the analytic tertiary invariant
$$\eta^{an}:S_{2k-2}\to  \frac{\C[[q]]}{E^\Gamma_{\C,2k}[[q]] + {}^N\Z [[q]]+\C}\ ,\quad   \eta^{an}([Z]):=[e_{ell}(N_0,Z)] \ ,$$
where we assume that $k\ge 2$. 
The main results of the present paper are the construction  of a topological
analog 
$$\eta^{top}:S_{2k-2}\to  \frac{\C[[q]]}{E^\Gamma_{\C,2k}[[q]]+ {}^N\Z [[q]]+\C}$$
 and the tertiary index Theorem \ref{dsadioasdlsadasdsadsadd}
\begin{equation}\label{fzifwefwfwefwfwfwef}
\eta^{an}=\eta^{top}\ .
\end{equation}
The construction of $\eta^{top}$ is quite involved. The details will be given in 
Section \ref{wqdqwodqwdqwdqwd}, culminating in Definition \ref{etatop}. The specialist will recognize that  on the topological side we try to perform the analogous constructions as in  the definition of $\eta^{an}$ on the analytic side.
The equality (\ref{fzifwefwfwefwfwfwef}) seems to be the first tertiary index theorem in the mathematical literature.
The basic principle of the  construction of tertiary invariants presented here also works in other situations. This will be demonstrated elsewhere.

The main idea  of the proof of the tertiary index theorem  is to relate both sides of (\ref{fzifwefwfwefwfwfwef})
to a third invariant, the $f$-invariant defined by Laures. The derivation of these relations is the content of Sections \ref{ufwiefwefwefwef} and \ref{csacoascascs}, while the definition of the $f$-invariant will be recalled in detail in Section \ref{uiewdwedwed}. For the pupose of the introduction
let us review some interesting homotopy theoretic aspects.
The key tool for computing the stable homotopy groups $\pi_*^s(S^0)$ is the 
Adams-Novikov spectral sequence 


\begin{equation}\label{ANSS} 
E_2^{s,t}=\Ext_{MU_*MU}^s(MU_*,\Omega^{t/2}MU_*)\Rightarrow\pi_{t-s}^s(S^0),
\end{equation}

cf. \cite{ravenel}, which defines a separated and exhaustive filtration
\begin{equation}\label{filtration} \pi_*^s(S^0)=F^0\supseteq F^1\supseteq\ldots\end{equation}
and homomorphisms 
\begin{equation}\label{uifuwiefwefewfwef}
 F^0/F^1\to E_2^{0,*},\quad 
 F^1/F^2\to E_2^{1,*}\ , \quad \mbox{and}\quad  F^2/F^3\to E_2^{2,*}.
\end{equation}
Here $MU_*$ denotes the bordism ring of stably almost complex 
manifolds. 
It is canonically a comodule for the Hopf algebroid $(MU_*,MU_*MU)$, and
the $\Ext$-group is  calculated in the abelian category of comodules.
The algebraic approximation $E_2^{s,*}$ to $\pi^s_*(S^0)$ is known completely only for 
$s\leq 2$ and represents the current edge of computational knowledge about $\pi_*^s(S^0)$,
c.f. for example \cite{k2local}. We have $E_2^{0,*}=E_2^{0,0}=\Z$, the groups $E_2^{1,t}$ are finite cyclic
with order given by denominators of Bernoulli numbers, and $E_2^{2,*}$ is very complicated but 
known explicitly by \cite{mrw}. A conceptual interpretation of $E_2^{2,*}$ in terms of congruences between
elliptic modular forms was only achieved recently \cite{mark2line} using the topological modular forms
of Goerss, Hopkins and Miller. 

Knowing  $E_2^{i,*}$  for $0\leq i\leq 2$ 
the natural next question  is which elements  are permanent cycles in
(\ref{ANSS}): $E_2^{0,0}$ is permanent for trivial reasons and detects $\pi^s_0(S^0)\simeq\Z$. Deciding
which elements of $E_2^{1,*}$ are permanent is tantamount to Adam's famous solution of the Hopf
invariant one problem and $E_2^{1,*}$ exactly detects the image of the $J$-homomorphism $im(J)\subseteq\pi_*^s(S^0)$. A lot is known about the permament cycles in  $E_2^{2,*}$, with recent progress due to  the solution of the  Kervaire-invariant one problem \cite{hhr}.

Via the Pontrjagin-Thom isomorphism a closed $n$-dimensional framed manifold $X$ represents a class $[X]\in \Omega^{fr}_n\cong \pi^s_n(S^0)$. It is an interesting question to descide by which element in $E_2^{*,n+*}$ it is detected under the homomorphisms (\ref{uifuwiefwefewfwef}). If
$[X]$ represents a non-trivial element in $F^0/F^1$, then $n=0$ and the corresponding element of $E^{0,0}_2$ can
easily be calculated by counting points.
If $[X]$ represents a non-trivial element in $F^1/F^2$, then 
we have $n=2k-1$ for some $k\ge 1$.
In this case one can determine the corresponding element in $E^{1,2k}_2$ by comparing the $e$-invariant of $[X]$ with the known $e$-invariants of the elements
of $E^{1,2k}_2$. In \cite{deninger} it is demonstrated that the intrinsic analytic
expression $e^{an}$ of the $e$-invariant can by used to effectively calculate the bordism classes of certain framed nil-manifolds and to show that they  account for all of $im(J)$ (up to a factor of $2$).

If $n=2k-2$ with $k\ge 2$, then $[X]\in F^2$ automatically. 
In principle, the element in $E^{2,2k}_2$ represented by $[X]$ can be determined by calculating  the $f$-invariant of Laures \cite{laures-cob}. The receipe
given in  \cite{laures-cob} roughly requires to represent $X$ as a corner of codimension two of an almost complex manifold (the precise statement is explained in Section \ref{uiewdwedwed}). The relation of the $f$-invariant with the tertiary index theory invariant $\eta^{an}$ provides a first step towards an intrinsic
formula as it only requires to represent $X$ as a boundary of a (stably) almost complex manifold.  An honest intrinsic formula for the $f$-invariant is still unknown.
Nevertheless, already this first step can simplify calculations.  This has been demonstrated nicely by the explicit examples  calculated in the thesis \cite{bod}, though 
the precise analytic situation in this reference is different from the one considered in the present paper and more special. The approach of  \cite{bod} is based on 
manifolds with corners and  boundary fibration structures and associated eta-forms.
Using adabatic limits one can relate, or even derive the formulas for  the $f$-invariants
in  \cite{bod} from our  $\eta^{an}$.  At the moment we are not able to add
any new explicitly calculable example to the list of  \cite{bod}. Note that the comparison of  two $f$-invariants given by  formal power series  representatives in $\C[[q]]$  still leads to a very complicated computational problem
in the quotient  $\frac{\C[[q]]}{E^\Gamma_{\C2,k}[[q]] + {}^N\Z [[q]]+\C}\ .$

\textit{While working on this project we profited from discussions with G. Laures and Ch. Bodecker.}

\section{Dirac operators and the $e$-invariant}\label{duwiqdqwd}

In this Section we define the secondary invariants $e^{top}$ and $e^{an}$ mentioned in the Introduction. This  analytic interpretation
of Adams' $e$-invariant is due to Atiyah-Patodi-Singer \cite{MR0397798}. The main purpose
of this section is to set up notation and to explain the basic idea of the derivation of a secondary invariant from a primary one. The same principles will we applied in a much more complicated situation in the construction of  $\eta^{an}$ and $\eta^{top}$.

If $M$ is a closed almost complex manifold, then for every choice of a hermitean metric on $TM$ and a metric connection $\nabla^{TM}$ preserving the almost complex structure on  $TM$  the integral 
\begin{equation}\label{qwdqd}
\int_M \Td(\nabla^{TM})\in \R
\end{equation} of the Todd form is an integer, where
$$ \Td(\nabla^{TM})=\det\frac{\frac{R^{TM}}{2\pi i}}{1-\ee^{-\frac{R^{TM}}{2\pi i}}}$$
and 
$R^{TM}$ denotes the curvature form of $\nabla^{TM}$.
This follows from the Atiyah-Singer  index theorem
$$\ind(\Dirac_M)=\int_M \Td(\nabla^{TM})\ ,$$
where $\Dirac_M$ is the $Spin^c$-Dirac operator associated to the $Spin^c$-structure naturally induced by the almost complex structure.

If the manifold has a boundary $N=\partial M$, then in general the integral 
(\ref{qwdqd}) is just a real number. By the Atiyah-Patodi-Singer index theorem the combination
\begin{equation}\label{qwdqd1}
\int_M \Td(\nabla^{TM})+\left[\eta(\Dirac_N)+\int_N  \tilde \Td(\nabla^{LC,L},\nabla^{TM})  \right]
\end{equation}
is an index and therefore an integer, where
$\eta(\Dirac_N)\in \R$ is the $\eta$-invariant of the
 $Spin^c$-Dirac operator $\Dirac_N$ and
$\tilde \Td(\nabla^{LC,L},\nabla^{TM})$
is the transgression form which we explain in the following.
 The $\eta$-invariant is a global spectral invariant of $\Dirac_N$
and depends on the choice of a $Spin^c$-connection on $N$.
The group $Spin^c(n)$ fits into a central extension
$$1\to U(1)\stackrel{c}{\to} Spin^c(n)\to SO(n)\to 1\ .$$
Furthermore, there exist a homomorphism $u:Spin^c\to U(1)$ such that the composition
$u\circ c:U(1)\to U(1)$ is the double covering.
Therefore, a $Spin^c$-connection is determined by the Levi-Civita connection $\nabla^{LC}$ of the Riemannian metric and the central part
$\nabla^{L^2}$, a connection on the line bundle canonically associated to the $Spin^c$-structure via the character $u$. 
We have the following diagram of classical groups
$$\xymatrix{&U(1)\ar[dr]^2\ar[d]&\\
U(n)\ar[r]\ar@/^3cm/[rr]_{\det}&Spin^c(2n)\ar[d]\ar[r]^u&U(1)\\
&SO(2n)&}$$
which shows the following:
\begin{enumerate}
\item An  almost complex structure and a hermitean metric on $TM$, i.e. an $U$-structure, 
induces naturally a $Spin^c$-structure.
\item In this case the line bundle $L^2\to M$ given by the $Spin^c$-structure is $L^2\cong \Lambda^{m}_\C T^*M$.
 \end{enumerate}
If the $Spin^c$-structure
comes from an almost complex structure, then a  connection on $TM$ which preserves the metric and the  almost complex structure induces a connection on $L$. Note that
$\nabla^{LC}$ in general does not preserve the almost complex structure and therefore does not induce a connection on $L^2$.

The transgression of the Todd form   in (\ref{qwdqd1})
has the following precise meaning. We split
$$\frac{x}{1-\ee^{-x}}=\ee^{\frac{x}{2}}\frac{x/2}{\sinh(x/2)}\ .$$
The second factor is an even power series and gives a characteristic form 
$$\hA(\nabla^{TM})=\det^{1/2}\left(\frac{\frac{R^{TM}}{4\pi}}{\sinh(\frac{R^{TM}}{4\pi})}\right)$$
of the real bundle $TM$. The first factor 
$$\ch(\nabla^L)=\ee^{\frac{R^{TM}}{4\pi i}}$$ represents the
Chern character of a formal square root of the canonical bundle
$L^2=\Lambda^{m}T^*M$, if $\nabla^{TM}$ preserves the almost complex structure and the hermitean metric.
In this way we can rewrite the Todd-form
as a characteristic form associated to a pair $(\nabla^{TM},\nabla^{L^2})$ of a real connection on $TM$ and a connection on $L^2$.  
A metric complex connection
$\nabla^{TM}$ naturally gives rise to such a pair
$(\nabla^{L^2},\nabla^{TM})$, and in this case we have
$$\Td(\nabla^{TM})=\ch(\nabla^{L})\wedge \hA(\nabla^{TM})\ .$$ A $Spin^c$-connection gives rise to another pair
$(\nabla^{LC},\nabla^{L^2})$, and in this case we write
$$\Td(\nabla^{LC,L})=\ch(\nabla^{L})\wedge  \hA(\nabla^{LC})$$
The transgression form
 $\tilde \Td(\nabla^{LC,L},\nabla^{TM})$ interpolates between these ends in the sense that
$$d\tilde \Td(\nabla^{LC,L},\nabla^{TM})=
\Td(\nabla^{LC,L})-\Td(\nabla^{TM})\ .$$

The upshot of this discussion is that the class
$$[\int_M \Td(\nabla^{TM})]\in \R/\Z$$
is equal to $$[\int_N  \tilde \Td(\nabla^{TM},\nabla^{LC,L})  -\eta(\Dirac_N)]$$
and therefore only depends on the boundary $N$ of $M$ as a geometric object.

Let us now assume that the boundary is framed, i.e. we have fixed an isomorphism $TN\cong N\times \R^{2m-1}$, where
$2m=\dim_\R M$. Adding the normal direction we get an induced framing  $TM_{|N}\cong N\times \R^{2m}$ and, using $\R^{2m}\cong \C^m$, a metric and an almost complex structure induced by the framing. We assume that the given almost complex structure and metric on $TM$ restrict to the ones induced by the framing over $N$.
Furthermore we assume that the metric complex connection $\nabla^{TM}$ restricts to the trivial one $\nabla^{triv}$ over $N$. Then
$ \tilde \Td(\nabla^{LC,L},\nabla^{TM})_{|N}=\tilde \Td(\nabla^{LC,L},\nabla^{triv})$ does not depend on the remaining choice of $\nabla^{TM}$ at all. We conclude that in this case, the classes appearing in  (\ref{qwdqd1})
\begin{eqnarray}
e^{top}(N)&:=&[\int_M \Td(\nabla^{TM})]\in \R/\Z\label{zuzu1}\\ e^{an}(N)&:=&[\int_N  \tilde \Td(\nabla^{triv},\nabla^{LC,L})  -\eta(\Dirac_N)]\in \R/\Z\label{zuzu2}
\end{eqnarray}
are equal, i.e.
\begin{equation}\label{anauiw}
e^{an}(N)=e^{top}(N)\ ,\end{equation} and
that they only depend on the framed manifold $N$.
From now on we omit the superscripts $top$ and $an$.

It is easy to see that
$e(N)$ is a framed bordism invariant. In fact, the intrinsic interpretation (\ref{zuzu2}) shows that
$e(N\sqcup N^\prime)=e(N)+e(N^\prime)$.
If $M$ is a framed bordism between 
$N$ and $N^\prime$, then we can choose
the trivial connection $\nabla^{TM}:=\nabla^{triv}$ and therefore by (\ref{zuzu1})
$$e(N)-e(N^\prime)=e(N\sqcup -N^\prime)=[\int_M \Td(\nabla^{TM})]=0\ .$$
The Todd class is stable, i.e. if we add a trivial bundle
$V\cong M\times \R^r$ to $TM$ and let $\nabla^{V}$ be the trivial connection, then
$$\Td(M)=\Td(M\oplus V)\ , \ \Td(\nabla^{TM})=\Td(\nabla^{TM\oplus V}) \ .$$
A stable framing or stable almost complex structure on $M$ is a framing or almost complex structure on $TM^s:=TM\oplus V$ for a suitable $r$. A stable almost complex structure still induces a $Spin^c$-structure, and the discussion above easily extends to the stable setting. In particular, we get a homomorphism
$e:\Omega^{fr}_*\to \R/\Z$ from the bordism group of stably framed manifolds.

By the Pontrjagin-Thom construction the   group $\Omega_*^{fr}$  
is isomorphic to the stable homotopy group $\pi_*^S$ of the sphere.
If a class $[f]\in \pi_n^S$ is represented by a differentiable map
$f:S^{m+n}\to S^m$, then for a regular point
$x\in S^m$  the preimage
$N:=f^{-1}(\{x\})\subset S^{m+n}$ is an $n$-manifold whose stable normal bundle is framed. This framing induces an equivalence class of stable framings of the tangent bundle, and the corresponding
$[N]\in \Omega^{fr}_n$ represents
the image of $ [f]$ under the Pontrjagin-Thom isomorphism
$$\pi_n^S\stackrel{\sim}{\to} \Omega_n^{fr}\ .$$

The $e$-invariant
$$e:\pi_*^S\cong \Omega_*^{fr}\to \R/\Z$$ has been introduced by Adams \cite{adams} and was identified with the analytic expression (\ref{zuzu2}) by Atiyah-Patodi-Singer \cite[Theorem 4.14]{MR0397798}.

\section{Modular Dirac operators and $\eta^{an}$}\label{sdvdsv}

In this Section we first  recall the construction of the complex elliptic genus \cite{hirzebruch} and introduce the necessary notation in order to write down the corresponding formal power series of Dirac operators and its spectral invariants.
Then we introduce the analytic tertiary invariant $\eta^{an}$ adopting an innocent simplifying assumption. In the more technical Section \ref{dusadasdasd} 
this assumption we be removed, and the analytic derivation of the properties of $\eta^{an}$ will be given.

We fix a number $4\le N\in \nat$ and a primitive root of unity $\zeta_N$.  We consider the group $$\Gamma:=\Gamma_1(N):=\{\left(\begin{array}{cc}a&b\\c&d\end{array}\right)|a,d\equiv 1(N)\ ,c\equiv 0(N)\}\subset SL(2,\Z)\ .$$
By $E^\Gamma_\C$ we denote the ring of modular forms for
$\Gamma$. Note that the group
$\Gamma$ acts on the upper half plane
$H=\{z\in\C\, | \, \Imm(z)>0\}$ by fractional linear transformations.
The quotient $\cM:=\Gamma\backslash H$ parameterizes elliptic curves with a distinguished point of order $N$. There is a universal elliptic curve
$u:\cE\to \cM$ with zero section
$e:\cM\to \cE$. The pull-back of the vertical bundle
$\bar \omega:=e^*Tu$ is a holomorphic line bundle which satisfies $\bar \omega^2=T^*\cM$ (Kodaira-Spencer).
Its lift $\omega$ to the upper half plane therefore is a $\Gamma$-equivariant square root of the canonical bundle
$T^*H$. A modular form of weight $k\in\Z$ for the group $\Gamma$ is a holomorphic section of  $\omega^k$ which is $\Gamma$-invariant and of moderate growth in the cusps. The ring $E^\Gamma_\C$ is non- negatively 
graded by the weight and of finite type, i.e. $\dim(E^\Gamma_{\C,k})<\infty$ for all $k\ge 0$. If one trivializes the bundle
$\omega^k$ by $(dz)^{k/2}$, then one identifies modular forms with functions on $H$. 
If we use the coordinate $q=\ee^{2\pi  i\tau}$, $\tau\in H$, then a modular form $\phi\in E_\C^\Gamma$ has a Fourier expansion
$\phi(q)=\sum_{n\ge 0} a_n q^n$. Following conventions
in topology, we will write $E^\Gamma_{\C,2k}$ for the space of
modular forms of weight $k$.

\begin{ddd}
We consider the ring $$\ZN:=\Z[\frac{1}{N},\zeta_N]$$ and call a modular form $\phi\in E^\Gamma_{\C,2k}$ of weight $k$ integral, if the coefficients in the expansion $\phi(q)=\sum_{n\ge 0} a_n q^n$  belong to $\ZN$. We let $E^\Gamma\subseteq E^\Gamma_\C$ denote the graded subring of integral modular forms.
\end{ddd}
 We consider the power series in $q$ and $x$, c.f. \cite[page 175]{hirzebruch}
$$Q_{y}(x)(q):=\frac{x}{1-e^{-x}} (1+ye^{-x})\prod_{n=1}^\infty\frac{1+yq^ne^{-x}}{1-q^ne^{-x}}\frac{1+y^{-1}q^ne^x}{1-q^ne^x}\ .$$
We further define
$$a(q):=Q_{-\zeta_N}(0)(q)^{-1}$$
and \begin{equation}\label{uiidqwdwqd}
\phi(x)(q):=a(q)Q_{-\zeta_N}(x)(q).
\end{equation}
Then the following is known from the classical theory of theta-functions:
\begin{lem}\label{zuidqwdqwd}
If we expand \begin{equation}\label{wdqwdqdqwd}\phi(x)(q)=\sum_{n\ge 0} \phi_{n}(q)x^n\end{equation}
then $\phi_n(q)$ is the $q$-expansion of a  modular form $\phi_n\in E^\Gamma_{\C,2n}$ of weight $n$. Moreover, $\phi_0=1$.
\end{lem}

Let now $M$ be an almost complex manifold of real dimension $2n$. If we choose a hermitean metric and a    connection $\nabla^{TM}$ preserving the almost complex structure and the metric then
we can define the element
$$\phi(\nabla^{TM}):=\det(\phi(\frac{R^{TM}}{2\pi i}))\in \Omega(M)\otimes E^\Gamma_\C\ .$$
More precisely, we write
$$\prod_{i=1}^n\phi(x_i)(q)=
 \sum_{n\ge 0}K_{n}(\sigma_1,\dots,\sigma_n)\psi_n(q)\ ,$$
where $K_n$ is homogeneous of total degree $n$ and
$\psi_n\in E^\Gamma_{\C,2n}$ is a homogeneous polynomial of total degree $n$ in the modular forms $\phi_k$ appearing in (\ref{wdqwdqdqwd}).  The 
$\sigma_i:=\sigma_i(x_1,\dots,x_n)$ denote the elementary symmetric functions.  In terms of the Chern forms $c_i(\nabla^{TM})$ we have
\begin{equation}\label{ufiewfwef}
\phi(\nabla^{TM})_{2k}=K_k(c_1(\nabla^{TM}),\dots, ,c_n(\nabla^{TM})) \psi_{k}\in \Omega^{2k}(M)\otimes E^\Gamma_{\C,2k}\ .
\end{equation}
We now replace the Todd form in (\ref{qwdqd})
by $\phi(\nabla^{TM})$ and get the modular form
\begin{equation}
\phi(M):=\int_M\phi(\nabla^{TM})\in E^\Gamma_{\C,2n}\ .
\end{equation}

It again follows from an index theorem that this modular form is integral:
 \begin{lem}
We have $$\phi(M)=\int_M\phi(\nabla^{TM})\in E^\Gamma_{ 2n}\ .$$
\end{lem}
\proof
 We use the following calculus of power series with coefficients in the semigroup of vector bundles on $M$.
For a complex vector bundle $V\to M$ we consider the power series
 $$\Lambda_tV:=\sum_{i=0}^{\dim V} \Lambda^iV t^i\ ,\quad S_tW:=\sum_{i=0}^\infty S^iV t^i\ ,$$
where $\Lambda^i$ (resp. $S^i$) denotes the $i^{th}$ exterior (resp.
symmetric) power. 
If the $x_i$ denote the formal Chern roots of $V$ \footnote{The precise meaning of formal Chern roots is the following. One forms the bundle $\pi:F(V)\to M$ of complete flags in $V$. The pull-back by $\pi$ induces an injection $\pi^*:H^*(M;\Z)\hookrightarrow H^*(F(V);\Z)$.
  The pull-back  $\pi^*V$   has  a canonical  decomposition  $\pi^*V\cong \oplus_{i=1}^{\dim(V)} L_i$ as a sum of line bundles, and $x_i:=c_1(L_i)\in H^2(F(V);\Z)$. The elementary symmetric functions in the Chern roots are the pull-backs of the Chern classes of $V$, i.e. $\sigma_i(x_1,\dots,x_n)=\pi^*c_i(V)$. To be precise, the following formulas have to be interpreted in
$H^*(F(V);\Q)$},
then we have
$$\ch\Lambda_t V=\prod_i (1+te^{x_i})\ ,\quad \ch S_tV=\prod_i (1-t\ee^{x_i})^{-1}\ .$$ Furthermore we have 
$\Td(V):=\prod_i \frac{x_i}{1-\ee^{-x_i}}$.
It follows that $$\prod_i Q_y(x_i)=\Td(V)\:\: \ch\left[ \Lambda_y V^* \:\:\prod_{n=1}^\infty \Lambda_{q^ny}V^*\:\: \Lambda_{q^ny^{-1}} V \:\:S_{q^n}(V+V^*)\right] \ .$$
We form the formal power series in $q$
\begin{equation}\label{zuidwqdqwdwqdqwd}
C(V)(q):=a(q)^{\dim(V)}\Lambda_{-\zeta_N}(V^*)\prod_{n=1}^\infty\Lambda_{-\zeta_Nq^n}(V^*)\Lambda_{-\zeta_N^{-1}q^n}(V)S_{q^n}(V\oplus V^*) 
\end{equation}
 with coefficients in the semigroup of vector bundles and $\ZN$, i.e.
\begin{equation}\label{euwdioewdwed}C(V)(q)=\sum_{n\ge 0} W_n c_n q^n\ ,\end{equation}
where $W_n\to M$ is some vector bundle on $M$ functorially derived from $V$ (i.e. a combination of alternating and symmetric powers), and $c_n\in \ZN$.
A metric and a compatible  connection on $V$ naturally induces a metric and a compatible connection on all the coefficient bundles $W_n$.  Taking the Chern forms
we get
the formal power series
$$\ch(\nabla^{C(V)(q)}):=\sum_{n\ge 0} \ch(\nabla^{W_n}) c_n q^n \ .$$
 
In view of the definition (\ref{uiidqwdwqd}) we see that
$$\phi(\nabla^{TM})(q)=\Td(\nabla^{TM})\wedge \ch(\nabla^{C(TM)(q)})=\sum_{n\ge 0}\Td(\nabla^{TM})\wedge \ch(\nabla^{W_n}) c_nq^n\ .$$
A hermitean vector bundle with a compatible connection $(W,\nabla^W)$ can be used to form the twisted Dirac operator $\Dirac_M\otimes W$. The formal power series
$$\Dirac_M\otimes C(V)(q):=\sum_{n\ge 0}c_nq^n\Dirac_M\otimes W_n $$
of twisted Dirac operators is the modular Dirac operator alerted to in the title.
The Atiyah-Singer index theorem gives
$$\ind(\Dirac_M\otimes W_n)=\int_M \Td(\nabla^{TM})\wedge \ch(\nabla^{W_n})\in \Z\ .$$
This implies that the expansion
$$\int_M \phi(\nabla^{TM})(q)=\sum_{n\ge 0}c_nq^n \ind(\Dirac_M\otimes W_n )$$
has coefficients in $\ZN$, and we conclude that
$$\phi(M)=\int_M \phi(\nabla^{TM})\in E^\Gamma_{2n}\ .$$
\hB
 By construction we have
$\phi(M_0\cup M_1)=\phi(M_0)+\phi(M_1)$.
For a product $M_0\times M_1$ we choose the product connection on
$\pr_0^*TM_0\oplus \pr_1^*TM_1$. Then we have
$$\phi(\nabla^{T(M_0\times M_1)})=\pr_0^*\phi(\nabla^{TM_0})
\wedge \pr_1^*\phi(\nabla^{TM_1})\ .$$
This implies that
$\phi(M_0\times M_1)=\phi(M_0)\phi(M_1)$.
 Finally, if $M$ is zero-bordant as a stably almost complex manifold, then $\phi(M)=0$ by Stokes' theorem.
We therefore obtain a homomorphism of graded rings
$\phi:MU_*\to E^\Gamma_{*}$.

\begin{ddd}\label{uiddqw}
The ring homomorphism
$\phi:MU_*\to E^\Gamma_*$ is called  the complex elliptic genus of level $N$.
\end{ddd}

Since $\Td(\nabla^{LC,L})$ is cohomologous to $\Td(\nabla^{TM})$ we can write
$$\phi(M)=\int_{M}\phi(\nabla^{TM})=\int_M \Td(\nabla^{LC,L})\wedge \ch(\nabla^{C(TM)})\ .$$

 Let us now assume that  $M$ has a  boundary $N$. We will choose the metric on $M$ with a product structure.  The expression
$\int_M \Td(\nabla^{LC,L})\wedge \ch(\nabla^{C(TM)})$ now gives an inhomogeneous element in $\oplus_{n\ge 0} E^\Gamma_{\C,2n}$.
In order to define a homogeneous element containing the term
$\Td(\nabla^{LC,L})$, which is important since we want to apply local index theory,
we first observe (see (\ref{ufiewfwef})) that
$$[\Td(\nabla^{TM})\wedge \ch(\nabla^{C(TM)})]_{2n}\in \Omega(M)^{2n}\otimes E^\Gamma_{\C,2n}\ .$$
Using Stoke's  theorem we write
\begin{eqnarray}
\int_{M}\Td(\nabla^{TM})\wedge \ch(\nabla^{C(TM)})&=&
\int_{M}\Td(\nabla^{LC,L})\wedge \ch(\nabla^{C(TM)})\nonumber\\&&+\int_{M}d\tilde \Td(\nabla^{TM},\nabla^{LC,L})\wedge \ch(\nabla^{C(TM)})\nonumber\\
&=&\int_{M}\Td(\nabla^{LC,L})\wedge \ch(\nabla^{C(TM)})\nonumber\\&&+
\int_{N}\tilde \Td(\nabla^{TM},\nabla^{LC,L})\wedge \ch(\nabla^{C(TM)})\nonumber\\&\in& E^\Gamma_{\C,2n}\label{uidqwdqwd}\ ,
\end{eqnarray}
where $\tilde \Td(\nabla^{TM},\nabla^{LC,L})$ is the transgression of the Todd form satisfying
$$ d\tilde \Td(\nabla^{TM},\nabla^{LC,L})=\Td(\nabla^{TM})-\Td(\nabla^{LC,L})\ .$$

We again apply the Atiyah-Patodi-Singer index theorem to the twisted operators $\Dirac_M\otimes W_n$:
The sum
$$\int_M \Td(\nabla^{LC,L})\wedge \ch(\nabla^{W_n})+\eta(\Dirac_N\otimes W_{n|N})$$ is an index and therefore an integer.
Let us write
\begin{equation}\label{fuwehiewfewf}\eta(\Dirac_N\otimes C(TM_{|N})(q)):=\sum_{n\ge 0} c_nq^n 
\eta(\Dirac_N\otimes W_{n|N})\in \C[[q]]\ .\end{equation}
Then we have
$$\int_M \Td(\nabla^{LC,L})\wedge \ch(\nabla^{C(TM)(q)})+\eta(\Dirac_N\otimes C(TM_{|N})(q))\in \ZN[[q]]\ .$$
Therefore the Atiyah-Patodi-Singer theorem implies that 
\begin{equation}\label{udwqdqwdqwd}\int_{N}\tilde \Td(\nabla^{TM},\nabla^{LC,L})\wedge \ch(\nabla^{C(TM)(q)})-\eta(\Dirac_N\otimes C(TM_{|N})(q))\in E^\Gamma_{\C,2n}[[q]]+\ZN[[q]]\ ,\end{equation}
where
$$E^\Gamma_{\C,2n}[[q]]\subseteq \C[[q]]$$ denotes the finite-dimensional subspace of $q$-expansions of elements of $E^\Gamma_{\C,2n}$.
If $V\to N$ is a trivial bundle with the trivial connection and
$C(V)(q)=\sum_{n\ge 0} c_nq^n W_n$, then $W_n$ is trivial and $\eta(\Dirac_N\otimes W_{n|N})=\dim(W_n)\eta(\Dirac_N)$. Because of our normalization (\ref{uiidqwdwqd}) we have
$$\sum_{n\ge 0}c_nq^n\dim(W_n)=1\ .$$
 We conclude that for trivial $V$
 \begin{equation}\label{duqwdwqdqwd}\eta(\Dirac_N\otimes C(V)(q))=\eta(\Dirac_N)\ .\end{equation}
Similarly, $$\int_{N}\tilde \Td(\nabla^{TN},\nabla^{LC,L})\wedge \ch(\nabla^{C(V)(q)})=
\int_{N}\tilde \Td(\nabla^{TN},\nabla^{LC,L})\ .$$
Hence we have
$$\int_N \tilde \Td(\nabla^{TN},\nabla^{LC,L})\wedge \ch(\nabla^{C(V)(q)})-\eta(\Dirac_N\otimes C(V)(q))\in \C\subset \C[[q]]\ .$$

If we assume  that $N$ is framed and  that the almost complex structure and the connection on $TM$ are  compatible with the framing, then
$$\int_M \Td(\nabla^{LC,L})\wedge \ch(\nabla^{C(V)(q)})\in (E^\Gamma_{\C,2n}[[q]]+\C)\cap \ZN[[q]]\ .$$

Let us now consider the $2n-1$-dimensional manifold $N$ with a stable almost complex structure as the primary object. After choosing a Riemannian metric and a $Spin^c$-connection we can define
$$\int_{N}\tilde \Td(\nabla^{TN^s},\nabla^{LC,L})\wedge \ch(\nabla^{C(TN^s)(q)})-\eta(\Dirac_N\otimes C(TN^s)(q))\in \C[[q]]\ ,$$
where $TN^s\cong TN\oplus (N\times \R^k)$ denotes a stabilization of $TN$ which carries the almost complex structure and a complex connection $\nabla^{TN}$. The class
$$[\int_{N}\tilde \Td(\nabla^{TN},\nabla^{LC,L})\wedge \ch(\nabla^{C(TN)(q)})-\eta(\Dirac_N\otimes C(TN^s)(q))]\in \C[[q]]/\C$$ 
is invariant under further stabilization, i.e., under replacing
$TN^s$ by $TN^s\oplus (N\times \C^l)$ (where the second summand has the trivial connection).

Now observe that the bordism groups $MU_*$ of stably almost complex manifolds are concentrated in even degrees. Therefore
$MU_{2n-1}=0$, and $N$ admits a zero bordism $M$ with a stable almost complex structure. The discussion above implies that 
\begin{equation}\label{sdasdasdasdq} 0=[\int_{N}\tilde \Td(\nabla^{TN},\nabla^{LC,L})\wedge \ch(\nabla^{C(TN)(q)})-\eta(\Dirac_N\otimes C(TN^s)(q))]\in \frac{\C[[q]]}{ E^\Gamma_{\C,2n}[[q]]+\ZN[[q]]+\C}\ .\end{equation}
From the point of view of the spectral theory on $N$, this fact
is completely mysterious.

This equation is the higher analog of  the relation
$$0=[\int_{TM}\Td(\nabla^{TM})]\in \R/\Z$$ in the even-dimensional case.
 If $N$ has a boundary, then the equality (\ref{sdasdasdasdq})
is no longer true in general, and this defect is the principal topic of the present paper.

We now introduce one of the main objects of our investigations,
namely an invariant $\eta^{an}(Z)$ of a framed manifold $Z$ of
positive even dimension. The construction of this invariant in
full generality is somewhat technical and is deferred to Section 
\ref{xsjxkasxasx}. The suspicious reader will have to skip
ahead to Section \ref{xsjxkasxasx} now since we will use $\eta^{an}(Z)$
in the following. For the time being, we content ourselves with giving
the construction in a special case which reveals all the essential 
features.

In the above situation, we now consider the case that $N$ has a boundary $Z:=\partial N$ such that $TN^s_{|Z}$ is framed, and the almost complex structure is compatible with this framing.
 Furthermore we assume that the Riemannian metric $g^N$ has a product structure near $Z$.
 For simplicity let us assume here that $\Dirac_Z$ is invertible. This assumption will be dropped later in the technical Section \ref{xsjxkasxasx} using the notion of a taming. The restrictions $W_{n|Z}$ are now trivialized so that $\Dirac_Z\otimes W_{n|Z}$ is invertible
for all $n\ge 0$. 
In this case, using global Atiyah-Patodi-Singer boundary conditions, 
we get a selfadjoint extension of $\Dirac_N\otimes W_n$
and we can define the $\eta$-invariant
$\eta(\Dirac_N\otimes W_n)\in \R$ and therefore
$$\eta(\Dirac_N\otimes C(TN^s)(q))\in \C[[q]]\ .$$
Using an extension of the Atiyah-Patodi-Singer index theorem to manifolds with corners \cite{math.DG/0201112} we will show the following theorem. 
\begin{theorem}\label{etaan}
In the above situation, the element
\begin{eqnarray*}
\eta^{an}(Z)&:=&[\int_{N}\tilde \Td(\nabla^{TN},\nabla^{LC,L})\wedge \ch(\nabla^{C(TN)(q)})-\eta((\Dirac_N\otimes C(TN^s)(q)))]\\&\in& \frac{\C[[q]]}{E^\Gamma_{\C,2n}[[q]]+\ZN[[q]]+\C}=:U^q_{2n}
\end{eqnarray*}
only depends on the framed bordism class of $Z$ and defines a homomorphism
$$\eta^{an}:\pi^S_{2n-2}=F^2\pi^S_{2n-2}\to \frac{\C[[q]]}{E^\Gamma_{\C,2n}[[q]]+\ZN[[q]]+\C}\ $$
with ker$(\eta^{an})\subseteq F^3\pi^S_{2n-2}+F^2\pi^S_{2n-2}[N^\infty]$, where for an abelian group $A$ we write as usual $A[N^\infty]:=\{ a\in A\, |\{\exists k\in\nat | N^ka=0\}\}$
\end{theorem}


\section{A topological invariant $\eta^{top}$ and the index theorem}\label{wqdqwodqwdqwdqwd}

In this Section we work in the 
stable homotopy category in order to define the invariant $\eta^{top}$ of framed cobordism.
The construction of $\eta^{top}$ in a certain sense models step-by step on the topological side the construction of $\eta^{an}$. This and the way of bringing in $\Q/\Z$-versions of the corresponding cohomology theories  seems to be just one example of a general principle. As mentioned earlier we will discuss another example involving Adams operations in a future paper. 
We were guided by our experience with differential cohomology theories  in which the original cohomology theory and its $\R/\Z$-version are nicely combined, and which gives a very suitable formalism for the construction of secondary invariants like the $e$-invariant, see \cite{bunke-2007}. A corresponding theory suitable in a similar way for tertiary invariants like $\eta^{an}$ and $\eta^{top}$ has yet to be developed. The principles of the construction of $\eta^{top}$ are one of the main contribution of the present paper.


Let $MU$ denote the spectrum which represents the complex bordism homology theory. It is a ring spectrum with a unit
$\epsilon:S\to MU$, where $S$ is the sphere spectrum which represents the framed bordism homology theory.
We define the spectrum $\overline {MU}$ as the cofiber in the fiber sequence
$$S\stackrel{\epsilon}{\to}MU\to \overline{MU}\ .$$

A stable homotopy class $\alpha\in \pi_m^S$, $m>0$, is  a homotopy class of maps of spectra
$\alpha:\Sigma^m S\to S$, where $\Sigma^m S$ is the $m$-fold suspension of the sphere spectrum.
It fits into the following diagram.
\begin{equation}\label{duqwidwqdwqdqd}
\xymatrix{&\Sigma^{-1}MU\ar[d]\\&\Sigma^{-1}\overline{MU}\ar[d]\\  \Sigma^mS\ar@{:>}[ur]^{\hat \alpha}\ar@{.>}[dr]\ar[r]^\alpha&S\ar[d]^\epsilon\\&MU}\ .
\end{equation}
 Since $\pi^S_m$ is finite and $MU_m$ is torsion free the dotted arrow $\epsilon\circ \alpha$  is zero-homotopic.
Hence we get a lift $\hat \alpha\in \overline{MU}_{m+1}$ which is well-defined up to
the image of
$MU_{m+1}\to \overline{MU}_{m+1}$.
Let us now assume that $m$ is even and positive.  Then
$MU_{m+1}=0$ so that $\hat \alpha$ is actually unique. Furthermore,  $\overline{MU}_{m+1}$
is a finite group  isomorphic to $\pi_m^S$. 

Since $\Q$ is a flat abelian group the
association
$X\mapsto \overline{MU}_{\Q,*}(X):=\overline{MU}_*(X)\otimes \Q$ is again a homology theory. 
We let $\overline{MU}_\Q$ denote a spectrum representing $\overline{MU}_{\Q,*}(\dots)$. We have a natural homotopy class of maps
$\overline{MU}\to \overline{MU}_\Q$ and define $
\overline{MU}_{\Q/\Z}$ as the cofiber in
$$\overline{MU}\to \overline{MU}_\Q\to \overline{MU}_{\Q/\Z}\ .$$

We now consider the diagram
\begin{equation}\label{duqwdqwdqwd}\xymatrix{&\Sigma^{-2}\overline{MU}_\Q\ar[d]\\&\Sigma^{-2}\overline{MU}_{\Q/\Z}\ar[d]\\  \Sigma^mS\ar@{:>}[ur]^{\tilde \alpha_{\Q/\Z}}\ar@{.>}[dr]\ar[r]^{\hat \alpha}&\Sigma^{-1}\overline{MU}\ar[d]\\&\Sigma^{-1}\overline{MU}_\Q}\ .\end{equation}
Since $\hat\alpha$ is a torsion element the dotted arrow is zero homotopic, and we can choose
a lift $\tilde\alpha_{\Q/\Z}\in \overline{MU}_{\Q/\Z,m+2}$.
This element is well-defined up to
 the image of 
$\sigma:\overline{MU}_{\Q,m+2}\to \overline{MU}_{\Q/\Z,m+2}$.

We now prepare to use Landweber's exact functor theorem. Let $c$ denote a cusp of the congruence
sub-group $\Gamma_1(N)$ other than the cusp $\infty$, its existence is guranteed by \cite[Proposition 1.34, (iv)]{shimura}.
Let $E^\Gamma_*\subseteq \tilde{E}^\Gamma_*$
denote the graded ring of modular forms for $\Gamma_1(N)$ which are holomorphic except possibly at the cusp $c$.
We use the $MU_*$-module structure on $E^\Gamma_*$ given by the elliptic genus $\phi:MU_*\to E^\Gamma_*$ (see \ref{uiddqw}) in order to define the functor
\begin{equation}\label{uiiwqduwqdwdwqd}
X\mapsto \tilde{E}^\Gamma_*(X):= MU_*(X)\otimes_{MU_*}\tilde{E}^\Gamma_*
\end{equation}
from spaces to graded rings. The ring $\tilde{E}^\Gamma_*$ is not flat over $MU_*$, but 
it is Landweber exact, \cite[Theorem 6]{MR1235295}. We use the ring $\ZN$ where $N$ is inverted 
and the ring $\tilde{E}^\Gamma_*$ involving also some meromorphic modular forms
in order to ensure this property.  Landweber exactness implies that $\tilde{E}^\Gamma_*(\dots)$ is a homology theory and is represented by a spectrum $\tilde{E}^\Gamma$.


The transformation
$\kappa: MU_*(X)\to \tilde{E}^\Gamma_*(X)$, $x\mapsto x\otimes 1$, is represented by a morphism of ring-spectra
$\kappa : MU\to \tilde{E}^\Gamma$. 
By construction, for every space $X$ there is a factorization of $\kappa$
$$ MU_*(X)\longrightarrow MU_*(X)\otimes_{MU_*} E_*^\Gamma\subseteq \tilde{E}^\Gamma_*(X)\ ,$$
a fact which we will refer to informally by saying that 
the values of $\kappa$ are holomorphic at all
cusps, including $c$.

We need yet another homology theory called Tate homology, we refer
the reader to \cite[Sections 2.5 and 2.6]{ahs} for more details.
The underlying group-valued functor is given by 
$$X\mapsto T_*(X):=K_*(X)\otimes_\Z \ZN[[q]]$$ (this is indeed a homology theory since $\ZN[[q]]$ is flat over $\Z$), where $K_*$ is complex $K$-homology. There is a natural transformation
$\nu: MU_*(X)\to T_*(X)$ which has the following geometric description.
If the continuous map $f:M\to X$ from a closed almost complex manifold $M$ represents the class $[f]\in MU_*(X)$,
then
$$\nu([f])=f_*([M]_K\cap C(TM))\ ,$$
where we consider the formal power series
$C(TM)$ (see (\ref{zuidwqdqwdwqdqwd})) as an element of
$K^0(M)\otimes \ZN[[q]]$, 
$[M]_K$ is the $K$-theory fundamental class of $M$
(induced by the $Spin^c$-structure determined by the almost complex structure), and
$$\cap:K_*(M)\otimes (K^0(M)\otimes \ZN[[q]])\to K_*(M)\otimes \ZN[[q]]=T_*(M)$$ is the $\cap$-product between $K$-homology and $K$-theory. 

As a multiplicative homology theory Tate homology  is derived via the Landweber exact functor theorem from the formal group law of the Tate elliptic curve over $\ZN[[q]]$. This formal group law is classified by the homomorphism
$\nu:MU_*\to T_*$ defined above in the case $X:=*$.

We let $T$ denote a spectrum representing the Tate homology, and we use the symbol  $\nu:MU\to T$ also to denote a map of spectra representing the above transformation.
We now construct a map $\gamma:\tilde E^\Gamma\to T$ 
such that 
$$\xymatrix{MU\ar[rr]^\nu\ar[dr]^\kappa&&T\\&\tilde E^\Gamma\ar[ur]^\gamma&}$$
commutes up to homotopy:
We will construct the corresponding natural transformation of homology theories.
Note that $T_*$ is Landweber exact over
$MU_*$ so that we have a natural isomorphism
$$MU_*(X)\otimes_{MU_*}T_*\stackrel{\sim}{\to}T_*(X)$$
induced by $\nu\otimes 1$. Therefore in view of (\ref{uiiwqduwqdwdwqd}), in order to define a natural transformation of homology theories $\gamma$, we must only define a ring homomorphism 
$\gamma:\tilde E^\Gamma_*\to T_*$ such that
$\gamma\circ \kappa=\nu:MU_*\to T_*$.
The map
$$\gamma:\tilde E^\Gamma_{2n}\to K_{2n}\otimes \ZN[[q]]\cong \ZN[[q]]$$ which associates to the modular form
$\phi\in \tilde E^\Gamma_{2n}$ its $q$-expansion $\phi(q)\in \ZN[[q]]$
(and which is zero in odd degrees) has this property. Note that by definition all elements of $\tilde E^\Gamma_{2n}$ are holomorphic at the cusp $\infty$.

The homology theories $\tilde E^\Gamma_*$ and $T_*$ are multiplicative.
We define the spectra $\bar{\tilde E}^\Gamma$ and $\bar T$ again as the cofibers of the units
$$S\to \tilde E^\Gamma\to \bar{\tilde E}^\Gamma\ ,\quad S\to T\to \bar T\ .$$
Furthermore, we consider spectra
$\bar{\tilde  E}^\Gamma_\Q$ and $\bar T_\Q$ representing homology theories
$$\bar{\tilde E}^\Gamma_{\Q,*}(X)=\bar{\tilde E}^\Gamma_*(X)\otimes_\Z\Q\ ,\quad \bar T_{\Q,*}(X)=\bar T_*(X)\otimes_\Z\Q$$
and define
$\bar{ \tilde E}^\Gamma_{\Q/\Z}$ and 
$\bar T_{\Q/\Z}$ as the cofibers
$$\bar{\tilde E}^\Gamma\to \bar{\tilde E}^\Gamma_\Q\to \bar{\tilde E}^\Gamma_{\Q/\Z}\ ,\quad \bar T\to \bar T_\Q\to \bar T_{\Q/\Z}\ .$$

We have the following diagram
\begin{equation}\label{diag2}
\xymatrix{\bar T\wedge K\ar[r]^{q}&\bar T_{\Q}\wedge K& && \\\bar{\tilde E}\wedge K\ar[r]\ar[u]^{\bar \gamma\wedge \id}&\bar{\tilde E}^\Gamma_{\Q}\wedge K\ar[u]^{\bar \gamma_\Q\wedge \id}& & & \\\overline{MU}\wedge MU\ar[r]\ar[u]^{\bar \kappa\wedge \theta}&
\overline{MU}_\Q\wedge MU\ar[r]\ar[u]^{\bar \kappa_\Q\wedge \theta}&\overline{MU}_{\Q/\Z}\wedge MU\ar[r] &\Sigma \overline{MU}\wedge MU &
\\&
\overline{MU}_{\Q}\ar[r]^\pi\ar[u]^{\id\wedge \epsilon}&\overline{MU}_{\Q/\Z}\ar[r]\ar[u]^{\id\wedge \epsilon}&\Sigma \overline{MU}\ar[u]\ar[u]^{\id\wedge \epsilon}&\\&&\Sigma^{m+2}S,\ar@{.>}@/^4.3cm/[uuuul]_\eta\ar@{.>}@/_-0.5cm/[uul]^{\bar \eta}\ar[u]^{\Sigma^2\tilde \alpha_{\Q/\Z}}\ar@{.>}[ur]^{\Sigma^2\hat  \alpha}&&}
\end{equation}
where $\theta:MU\to K$ is the complex orientation of $K$-theory.

Let us explain the construction of the maps $\bar \kappa_\Q$ and $\bar \gamma_\Q$.   
First of all, $\kappa:MU\to \tilde E^\Gamma$ fits into
\begin{equation}\label{udiqwdqwdqwdwqd}\xymatrix{S\ar[r]\ar@{=}[d]&MU\ar[d]^\kappa\ar[r]&\overline{MU}\ar@{.>}[d]^{\bar \kappa}\ar[r]^\delta&\Sigma S\ar@{=}[d]\\ S\ar[r]&\tilde E^\Gamma\ar[r]&\bar{\tilde E}^\Gamma\ar[r]^{\delta'}&\Sigma S.}\end{equation}
The stable homotopy category is triangulated, and the horizontal lines are distinguished triangles. It follows from the general properties of a triangulated category that a map $\bar \kappa$ which fills this diagram exists. It is unique up to homotopy as we now show:
Assume $\bar\kappa'$ is a second lift
and consider $\nu:=\bar\kappa - \bar\kappa'$. Then there exists  an $\alpha:\Sigma S\to\bar{\tilde E}^\Gamma$
such that $\nu=\alpha\circ \delta$.
Since $\tilde E^\Gamma_1=0$ ($\tilde E^\Gamma$ is even) the
  canonical map $[\Sigma S,\bar{\tilde E}^\Gamma]\to
[\Sigma S,\Sigma S]\cong\Z$ is bijective. We write $n:=\delta'\circ \alpha\in\Z$. Since the right square in (\ref{udiqwdqwdqwdwqd}) commutes we get
$0=\delta'\circ\nu=\delta'\circ \alpha\circ \delta=n\delta$. We claim that this implies $n=0$.
If so, we see that $\alpha$ factors through some $\Sigma S\to \tilde E^\Gamma$, hence
$\alpha=0$ (since $\tilde E^\Gamma_1=0$) and $\nu=0$, as desired.

We show by contradiction that $n=0$. Let us assume that $n\not=0$. We first observe that for all $i\neq 0,1$ we have an exact 
sequence
\[ 0\to \tilde E^\Gamma_{i}\to \bar{\tilde E}^\Gamma_i\stackrel{\delta}{\to}S_{i-1}\to 0\]
since $\tilde E^\Gamma_*$ is torsion-free, and $S_k$ is finite for $k\ge 1$. On the other hand there exists $i\ge 2$ 
and an element $z\in S_{i-1}$ such that $nz\not=0$, in fact, such an 
element can be found in the image of the $J$-homomorphism, c.f.
\cite[Theorem 1.1.13]{ravenel}.
Let $\hat z\in \bar{\tilde E}^\Gamma_i$ be a preimage.
Then $0\not=nz=n\delta (\hat z) =0$ is the desired contradiction.

The construction of $\bar \gamma$ and $\bar \gamma_\Q$ is analogous.
Let us now explain the construction of the map $\bar \eta$. 
We have $\alpha\in F^2\pi^S_m$. This means that 
the lift $\hat \alpha\in \overline{MU}_{m+1}$ belongs to the kernel of the map
$$\overline{MU}_{m+1}\stackrel{\id\wedge \epsilon}{\to} (\overline{MU}\wedge MU)_{m+1}\ ,$$
or equivalently, that is admits a further lift $\tilde \alpha$ in the Adams resolution (\ref{ufwiefewfwef}) below.
Hence there exists a lift $\bar \eta\in (\overline{MU}_\Q\wedge MU)_{m+2}$ which is unique up to the image of
$(\overline{MU}\wedge MU)_{m+2}\to  (\overline{MU}_\Q\wedge MU)_{m+2}$.
If we fix the choice of $\tilde \alpha_{\Q/\Z}$, then the composition $$\eta:=(\bar \gamma_\Q\wedge \id)\circ (\bar\kappa_\Q\wedge \theta)\circ\bar \eta\in (\bar T_\Q\wedge K)_{m+2}$$ is well-defined up to elements in the image of
$$(\overline{MU}\wedge MU)_{m+2}\to (\bar {\tilde E}^\Gamma\wedge K)_{m+2}\to (\bar T_\Q\wedge K)_{m+2}\ .$$
When we incorporate the  indeterminacy of $\tilde \alpha_{\Q/\Z}$, then the class

%
%
%
%

\begin{equation}\label{asdasdasd}\hat \eta(\alpha)\in \frac{(\bar T_\Q\wedge K)_{m+2}}{q\circ (\bar \gamma\wedge \id) \circ (\bar\kappa\wedge\theta)(\overline{MU} \wedge MU)_{m+2}+ (\bar\gamma_\Q\wedge \id)\circ(\bar\kappa_\Q\wedge\theta)\circ(\id\wedge\epsilon)(\overline{MU}_{\Q,m+2})}\end{equation}
represented by $\eta$
is well-defined, i.e. it depends only on $\alpha\in\pi^S_m$.

We now calculate a suitable quotient of the group on the right-hand side of (\ref{asdasdasd}).
First of all $\bar T_{\Q,*}$ is concentrated in even degrees and we  have $$\bar T_{\Q,0}\cong \frac{{}^N\Z[[q]]\otimes\Q}{\Q(\zeta_N)}\ ,\quad \quad\quad\bar T_{\Q,2m}\cong {}^N\Z[[q]]\otimes\Q \ ,m\not=0$$
This gives
$$(\bar T_\Q\wedge K)_{m+2}\cong \frac{ {}^N\Z[[q_0]]\otimes\Q}{\Q(\zeta_N)}\oplus \bigoplus_{2s+2r=m+2\ ,s\not=0}{}^N\Z[[q_s]]\otimes\Q\ .$$ 
By \cite[Sec. 2.3]{laures-expansion}
the image of 
$q\circ (\bar \gamma\wedge \id)\circ  (\bar\kappa\wedge \theta):(\overline{MU}\wedge MU)_{m+2}\to (\bar T_\Q\wedge K)_{m+2}$ is contained in the  subgroup
$$\frac{\ZN[[q_0]]}{\ZN}\oplus \bigoplus_{2s+2r=m+2, s\not=0} \ZN[[q_s]]  \ .$$
Finally, 
$(\bar\gamma_\Q\wedge \id)\circ(\bar\kappa_\Q\wedge\theta)\circ(\id\wedge\epsilon)(\overline{MU}_{\Q,m+2})$
is contained in the subspace of $q_0$-expansions 
$E^\Gamma_{\Q,m+2}[[q_0]]$ of rational modular forms of weight $m+2$, again using that $\kappa$ takes holomorphic values.
Therefore we have constructed a well-defined
invariant
$$\hat \eta^{top}(\alpha)\in \frac{\frac{{}^N\Z[[q_0]]\otimes \Q}{\Q(\zeta_N)}\oplus \bigoplus_{2s+2r=m+2, s\not=0} {}^N\Z[[q_s]]\otimes \Q }{\frac{{}^N\Z[[q_0]]}{\ZN}\oplus \bigoplus_{2s+2r=m+2,s\not=0} {}^N\Z[[q_s]] +E^\Gamma_{\Q,m+2}[[q_0]]}\ .$$
The natural map
$E^\Gamma_{\Q,m+2} \to E^\Gamma_{\C,m+2}=E^\Gamma_{\Q,m+2}\otimes_\Q\C$ and the identification of all $q_s$ with a single variable $q$ induce
a natural map
$$  \frac{\frac{{}^N\Z[[q_0]]\otimes \Q}{\Q(\zeta_N)}\oplus \bigoplus_{2s+2r=m+2, s\not=0} {}^N\Z[[q_s]]\otimes \Q }{\frac{{}^N\Z[[q_0]]}{\ZN}\oplus \bigoplus_{2s+2r=m+2,s\not=0} {}^N\Z[[q_s]] +E^\Gamma_{\Q,m+2}[[q_0]]}
\to \frac{\C[[q]]}{\ZN[[q]]+E^\Gamma_{\C,m+2}[[q]]+\C}=U^q_{m+2}  $$
to the target of $\eta^{an}$.
\begin{ddd}\label{etatop}
For $m>0$ even, we let
$$\eta^{top}:\pi^S_{m}\to  \frac{\C[[q]]}{E^\Gamma_{\C,m+2}[[q]]+\ZN[[q]]+\C}=U^p_{m+2}$$
be the homomorphism
induced by $-\hat \eta^{top}$ ({\em sic !}) such that
$\eta^{top}(\alpha)\in U^q_{m+2}$ is the class represented by $-\hat \eta^{top}(\alpha).$
\end{ddd}

We can now state our index theorem:
\begin{theorem}\label{dsadioasdlsadasdsadsadd}
For even $m>0$ we have an equality of homomorphisms
$$\eta^{an}=\eta^{top}:\pi_{m}^S=F^2\pi^S_m\to  \frac{\C[[q]]}{E^\Gamma_{\C,m+2}[[q]]+\ZN[[q]]+\C}=U^q_{m+2} $$
with kernel contained in $F^3\pi^S_m+F^2\pi^S_m[N^\infty]$. 
\end{theorem}
This result will be proven in Section \ref{dusadasdasd} as Theorem \ref{dwqdqwdwqdwd}.

%
%
%
%
%

\section{The $f$-invariant}\label{uiewdwedwed}

The tertiary index Theorem \ref{dsadioasdlsadasdsadsadd} stating that $\eta^{an}=\eta^{top}$ 
will be proved  by relating the quantities on both sides of this equality with  the $f$-invariant of Laures, see  Definition \ref{deff}.
We will recall in detail the geometric as well as the homotopy theoretic description of the $f$-invariant given in \cite{laures-cob}. Both pictures will be needed in the two subsequent sections.
\\
The $f$-invariant takes values in a target which differs from the target of $\eta^{an}$ and $\eta^{top}$. The relation between these quantities will be obtained in several steps.
In the present section we begin with a step-by step reinterpretation of the 
$f$-invariant in a sequence of targets tending to the one of $\eta^{an}$ and $\eta^{top}$,
a  process which will be completed in the following two sections.  In this way we derive a  sequence of invariants which essentially contain the same information as the original $f$-invariant and will therefore be denoted by various variants of the symbol $f$ with decorations added\footnote{We apologize for introducing so much notation, but we want to avoid to use the same symbol for different objects.}.

Let us recall the construction of the 
canonical $MU$-based Adams resolution of 
the sphere spectrum $S$, c.f. \cite[Chapter 2,2]{ravenel}, i.e.
the following diagram.

\begin{equation}\label{ufwiefewfwef}\xymatrix{&\vdots&\vdots&\\&\Sigma^{-1}\overline{MU}\wedge\Sigma^{-1}\overline{MU}\wedge\Sigma^{-1}\overline{MU}\ar[r]\ar[d]&\Sigma^{-1}\overline{MU}\wedge \Sigma^{-1}\overline{MU}\wedge\Sigma^{-1}\overline{MU}\wedge MU \\&\Sigma^{-1}\overline{MU}\wedge\Sigma^{-1}\overline{MU}\ar[r]^{\id\wedge \id\wedge \epsilon}\ar[d]&\Sigma^{-1}\overline{MU}\wedge\Sigma^{-1}\overline{MU}\wedge MU\ar@{.>}[lu] \\&\Sigma^{-1}\overline{MU}\ar[r]^{\id\wedge \epsilon}\ar[d]&\Sigma^{-1}\overline{MU}\wedge MU\ar@{.>}[lu]_{\delta}\\\Sigma^m S\ar[r]^\alpha\ar@{-->}[ur]^{\hat \alpha}\ar@{-->}[uur]^{\tilde \alpha}&S\ar[r]^\epsilon&MU\ar@{.>}[ul]}\ .\end{equation}
The horizontal arrows are induced by the unit $\epsilon:S\to MU$, and the triangles are fiber sequences.
It follows from the construction of the Adams-Novikov spectral sequence that a class $\alpha:\Sigma^mS \to S$ belongs, for example, to $F^2\pi_m^S$, if 
and only if it admits a lift
$$\tilde \alpha:\Sigma^m S\to \Sigma^{-1}\overline{MU}\wedge\Sigma^{-1}\overline{MU}$$
as indicated (a similar assertion holds true for all steps of the filtration).
We now assume that $m>0$ is even which implies that $\alpha\in F^2\pi_m^S$.
We have already seen in (\ref{duqwidwqdwqdqd})
that the first lift $\hat \alpha$ is unique up to homotopy.
Therefore the lift $\tilde \alpha$ is determined up to the image of
$\delta:(\overline{MU}\wedge MU)_{m+2}\to (\overline{MU}\wedge \overline{MU})_{m+2}$.
 The composition (in order to simplify the notation we shift by two)
\begin{equation}\label{composition}
\Sigma^{m+2} S\stackrel{\tilde \alpha}{\to }\overline{MU}\wedge \overline{MU}\stackrel{\bar\kappa\wedge \bar \kappa}{\to} \bar {\tilde E}^\Gamma\wedge \bar {\tilde E}^\Gamma\to 
 \bar {\tilde E}_\Q^\Gamma\wedge \bar {\tilde E}_\Q^\Gamma
\end{equation}
determines a class in
$$(\bar {\tilde E}_\Q^\Gamma\wedge \bar {\tilde E}_\Q^\Gamma)_{m+2}=\frac{(\tilde{E}_\Q^\Gamma\otimes \tilde{E}^\Gamma_\Q)_{m+2}}{\tilde{E}_{\Q,m+2}^\Gamma\otimes \Q + \Q\otimes \tilde{E}_{\Q,m+2}^\Gamma}. $$
It was shown in \cite[Theorem 2.3.1]{laures-expansion}, that if $\tilde \alpha$ is in the image of $\delta$, then
it gives rise to a class in 
$$\tilde{E}^\Gamma_{m+2}\tilde{E}^\Gamma+ \tilde{E}_{\Q,m+2}^\Gamma\otimes \Q + \Q\otimes \tilde{E}_{\Q,m+2}^\Gamma\subseteq (\tilde{E}_\Q^\Gamma\otimes \tilde{E}_\Q^\Gamma)_{m+2} $$
(more precisely, $\tilde{E}^\Gamma_{m+2}\tilde{E}^\Gamma$ denotes image of this group in $(\tilde{E}_\Q^\Gamma\otimes \tilde{E}^\Gamma_\Q)_{m+2}$ under the natural map
$$\tilde{E}^\Gamma_*\tilde{E}^\Gamma\to \tilde{E}^\Gamma_*\tilde{E}^\Gamma\otimes\Q\cong \tilde{E}^\Gamma_{\Q,*}\otimes_\Q \tilde{E}^\Gamma_{\Q,*}\ ).$$
We have thus defined a map sending $\alpha$ to the composition
in (\ref{composition})

\begin{equation}\label{dwqgzduqwd}
f_\Q:F^2\pi^S_{m}\to  \frac{(\tilde{E}_\Q^\Gamma\otimes \tilde{E}^\Gamma_\Q)_{m+2}}{\tilde{E}^\Gamma_{m+2}\tilde{E}^\Gamma+\tilde{E}_{\Q,m+2}^\Gamma\otimes \Q + \Q\otimes \tilde{E}_{\Q,m+2}^\Gamma} =: V_{\Q,m+2}.
\end{equation}

This version of the $f$-invariant is already a derived one.
The universal $f$-invariant is given by the natural map, well-known
to be injective,
$$f_{univ}: F^2\pi^S_{m}/F^3\pi^S_{m}\hookrightarrow E^{2,m+2}_{2,MU}=\Ext^{2,m+2}_{MU_*MU}(MU_*,MU_*)\ ,$$
where the target is a component of the $E_2$-term of the $MU$-based Adams
spectral sequence (\ref{ANSS}). Since $\kappa:MU\to \tilde{E}^\Gamma$ is Landweber exact of height two, the induced map
$$\kappa:\Ext^{2,m+2}_{MU_*MU}(MU_*,MU_*)\to 
\Ext^{2,m+2}_{\tilde{E}^\Gamma_*\tilde{E}^\Gamma}(\tilde{E}^\Gamma_*,\tilde{E}^\Gamma_*)
$$
is injective after inverting $N$. Furthermore, there is an injective map
$$\iota:
\Ext^{2,m+2}_{\tilde{E}^\Gamma_*\tilde{E}^\Gamma}(\tilde{E}^\Gamma_*,\tilde{E}^\Gamma_*)\to
V_{\Q,m+2}.
$$
\textcolor{black}{This result is due to Laures and has been generalized by Hovey
to arbitrary chromatic level. See \cite[Section 2.3]{hornbostelnaumann} for a detailed 
account using the present set-up.}
The relation between the $f$-invariant and the
universal $f$-invariant is now given by
$$f_\Q=\iota\circ \kappa \circ f_{univ}\ .$$
We conclude that $f_\Q$ factors over the quotient $F^2\pi^S_{m}\to F^2\pi^S_{m}/F^3\pi^S_{m}$, and since $\iota\circ \kappa$ is injective after inverting $N$, $f_\Q$ induces an injection

$$\xymatrix{ f_\Q:(F^2\pi^S_{m}/F^3\pi^S_{m})\otimes_\Z\Z[\frac{1}{N}]\ar@{^{(}->}[r] &  V_{\Q,m+2}.}$$


The theory developed in \cite{laures-cob} attaches a geometric meaning to the choice of $\tilde \alpha$. If we represent $\alpha$
by a framed $m$-manifold $Z$, then the choice of $\tilde \alpha$ corresponds to the choice of the following data (here $TZ$, $TY$, etc. denote representatives of the stable tangent bundle) which exist according to \cite{laures-cob}:
\begin{enumerate}
\item a decomposition $TZ\cong T^0Z\oplus T^1Z$ of framed bundles 
\item compact manifolds $Y_0,Y_1$ with boundary $\partial Y_0\cong Z\cong -\partial Y_1$.
\item decompositions $TY_i\cong T^0Y_i\oplus T^1Y_i$ together with  complex structures on
$T^iY_i$ and framings on $T^{1-i}Y_i$ such that:
\item The inclusion $Z\hookrightarrow Y_i$ identifies
$(T^1Y_0)_{|Z}\cong T^1Z$ and $(T^0Y_1)_{|Z}\cong T^0Z$  as framed bundles, and
$(T^0Y_0)_{|Z}\cong T^0Z $ and $(T^1Y_1)_{|Z}\cong T^1Z $ as complex bundles.
\item  a manifold with corners $X$ such that $\partial_0X\cong Y_0$ and $\partial_1X\cong  Y_1$.
\item 
 a decomposition $TX\cong T^0X\oplus T^1X$ of complex bundles such that:
\item The inclusions $Y_i\hookrightarrow X$ identify  $T^0X_{|Y_0}\cong T^0Y_0$,  $T^1X_{|Y_1}\cong T^1Y_1$,
$T^1X_{|Y_0}\cong T^1Y_0 $ and $T^0X_{|Y_1}\cong T^0Y_1 $ as complex bundles.
\end{enumerate}
These data refine $Z$ into a representative of a class
$$[Z]\in \Omega^{(U,fr)^2}_{n+2}$$ in the language of \cite{laures-cob}.
Let us call this collection of data a $<2>$-manifold which extends the framed manifold $Z$.
The collection of 1.- 3. (i.e. forgetting $X$ and related structure) will be called
a $\partial <2>$-manifold which extends $Z$. Finally, $X$ will then be called
a $<2>$-manifold which extends
the $\partial <2>$-manifold data.

We choose hermitean metrics on $T^iX$ and  metric connections $\nabla^{T^iX}$ which preserve the complex structures
and coincide with the trivial connection induced by the framing when restricted to
$Y_{1-i}$. Recall the definition of $C(V)(q)$ in (\ref{zuidwqdqwdwqdqwd}).
We define
$$\hat F(X):=\int_X  \Td(\nabla^{TX}) \wedge \ch(\nabla^{C(T^0X)(p)})\wedge \ch(\nabla^{C(T^1X)(q)})\in \C[[p,q]] \ .$$
A priori, this is an element in
$\C[[p,q]]$, but because of Lemma \ref{zuidqwdqwd}
 we actually have (recall that $\dim(X)=m+2$) 
\begin{equation}\label{alpha}
\hat F(X)\in (E_\C^\Gamma\otimes_\C E_\C^\Gamma)_{m+2}[[p,q]]\subseteq\C[[p,q]]\ .
\end{equation}
   
 We define
$$V_{m+2}:=\frac{(\tilde{E}_\C^\Gamma\otimes \tilde{E}^\Gamma_\C)_{m+2}}{\tilde{E}^\Gamma_{m+2}\tilde{E}^\Gamma+\tilde{E}_{\C,m+2}^\Gamma\otimes \C + \C\otimes \tilde{E}_{\C,m+2}^\Gamma}
$$
and let $F(X)\in V_{m+2}$ be the class represented by $\hat F(X)$.
It is shown in \cite{laures-cob} that the class
$F(X)$ is the image of the $f$-invariant $f_\Q(Z)$ of the corner
$Z$ under the inclusion
$
V_{\Q,m+2}
\hookrightarrow V_{m+2}$.
It thus only depends on the framed bordism class of $Z$. 


We now consider the quotient
$$W_{\Q,m+2}:=\frac{\Q(\zeta_N) [[p,q]]}{\ZN[[p,q]]+ \tilde{E}^\Gamma_{\Q,m+2}[[q]]+ \tilde{E}^\Gamma_{\Q,m+2}[[p]]+\Q(\zeta_N)}\ .$$
Since the $p,q$-expansion maps (c.f. \cite[Section 2.3]{laures-expansion})
$\tilde{E}^\Gamma_{m+2}\tilde{E}^\Gamma$ to $\ZN[[p,q]]$ it induces a natural map
$$i^\Q:V_{\Q,m+2}\to W_{\Q,m+2}\ .$$ 

\begin{lem}
The composition 
$$\xymatrix{ i^\Q\circ f_\Q:(F^2\pi^S_m/F^3\pi^S_m)\otimes_\Z\Z[\frac{1}{N}] \ar@{^{(}->}[r] & W_{\Q,m+2}}$$ is injective.
\end{lem}
\proof
This proof is based on \cite[Lemma 3.2.2]{laures-expansion}.
We consider $\alpha\in F^2\pi^S_m/F^3\pi^S_m$ and assume that $i^\Q(f_\Q(\alpha))=0$.
Note that
 $$E^{2,m+2}_{2,\tilde{E}^\Gamma}:=\Ext^{2,m+2}_{\tilde{E}^\Gamma_*\tilde{E}^\Gamma}(\tilde{E}^\Gamma_*,\tilde{E}^\Gamma_*)$$ is a component of the $E_2$-term of the $\tilde{E}^\Gamma$
-based Adams-Novikov spectral sequence.
For $\alpha\in F^2\pi_m^S$ we have
$\kappa(f_{univ}(\alpha))\in E^{2,m+2}_{2,\tilde{E}^\Gamma}$. Let $\Phi\in (\tilde{E}^\Gamma_\Q\otimes \tilde{E}^\Gamma_\Q)_{m+2}$ be a representative of the image of  this cycle under $\iota$. By assumption
there are $u,v\in \tilde{E}^\Gamma_{\Q,m+2}$, $c\in \Q(\zeta_N)$ and
$z\in \ZN[[p,q]]$ such that
$\Phi(p,q)=z(p,q)+u(p)+v(q)+c$.
Let us write
$\Phi(p,q)=\sum_{i,j\ge 0} \Phi_{ij}p^iq^j$, $z(p,q)=\sum_{i,j\ge 0} z_{ij}p^iq^j$,
and 
$u(p)=\sum_{i\ge 0} u_i p^i$. 
Then, setting $p=0$ above, we conclude that
$$\sum_{j\ge 0} \Phi_{0j} q^j=\sum_{j\ge 0} z_{0,j} q^j+
v(q)+u_0+c\in \tilde{E}^\Gamma_{\Q,m+2}[[q]]+\Q(\zeta_N)+\ZN[[q]]\ .$$
By \cite[Lemma 3.2.2, (iv)$\Rightarrow$(ii)]{laures-expansion} we have
$$\Phi(p,q)\in \tilde{E}^\Gamma_{\Q,m+2}[[p]]+\tilde{E}^\Gamma_{\Q,m+2}[[q]]+
\ZN[[p,q]]\ ,$$ and hence that
$\iota(\kappa(f_{univ}(\alpha)))=0$.
From the injectivity of $\iota\circ \kappa\circ f_{univ}$
we conclude that $\alpha=0$.
\hB
Let us finally define
\begin{equation}\label{duidqwdqwd}W_{m+2}:=\frac{\C [[p,q]]}{\ZN[[p,q]]+ \tilde{E}^\Gamma_{\C,m+2}[[q]]+ \tilde{E}^\Gamma_{\C,m+2}[[p]]+\C}\ \end{equation}
and consider the obvious injection
$j:W_{\Q,m+2}\to W_{m+2}$
and the natural map
$i:V_{m+2}\to W_{m+2}$ given by the $(p,q)$-expansion.
Then
$i(F(X))=j(i^\Q(f_\Q(\alpha)))$.
The upshot of this discussion is the commutative diagram
\begin{equation}\label{belief}\xymatrix{F^2\pi_m^S/F^3\pi_m^S\ar@/^2.5cm/[rrdd]^F\ar@/_3cm/[rrddd]_{\tilde{f}}\ar[rdd]^{f_\Q}\ar[r]^{f_{univ}}&E^{2,m+2}_{2,MU}\ar[d]^{\kappa}&\\&E^{2,m+2}_{2,\tilde{E}^\Gamma}\ar[d]^{\iota}&\\
&V_{\Q,m+2}\ar[d]^{i^\Q}\ar[r]&V_{m+2}\ar[d]^i\\
&W_{\Q,m+2}\ar[r]^j&W_{m+2},}\end{equation}
where $\tilde{f}:=j\circ i^\Q\circ f_\Q$ is injective after inverting $N$.\\
We need to refine this construction slightly. Define
\[ \tilde{W}_{m+2}:=\frac{\C[[p,q]]}{\ZN[[p,q]]+E_{\C,m+2}^\Gamma [[q]]+E_{\C,m+2}^\Gamma [[p]]+\C}.\]
By the definition of $W_{m+2}$ in (\ref{duidqwdqwd}) there is a canonical surjection
$\tilde{\pi}:\tilde{W}_{m+2}\to W_{m+2}$. According to (\ref{alpha}), the class $i\circ F(X)$ has a holomorphic representative, and even better,  the diagramm 
$$\xymatrix{ &\frac{F^2\pi^S_m}{F^3\pi^S_m}\ar[dd]^{f_\Q}\ar@/^4cm/@{.>}[dddd]^{\tilde f}\ar[dl]^{[\alpha]\mapsto [\tilde \alpha]}\\\frac{(\overline{MU}_\Q\wedge \overline{MU}_\Q)_{m+2}}{\delta(\overline{MU}\wedge MU)_{m+2}}\ar[d]^{\cong} &\\\frac{(MU_{\Q}\otimes MU_{\Q})_{m+2}}{MU_{m+2}MU+MU_{m+2}\otimes \Q+\Q\otimes MU_{m+2}}\ar[r]^{\kappa\otimes \kappa}\ar[d]&\frac{(\tilde E^\Gamma_\Q\otimes \tilde E^\Gamma_\Q)_{m+2}}{\tilde E^\Gamma_{m+2}\tilde E^\Gamma+\tilde E^\Gamma_{\Q,m+2}\otimes \Q+\Q\otimes \tilde E^\Gamma_{\Q,m+2}}\ar[d]^{j\circ i^\Q}\\
\frac{\C[[p,q]]}{\ZN[[p,q]]+  E^\Gamma_{\C,m+2}[[p]]+  E^\Gamma_{\C,m+2}[[q]]+\C}\ar[r]\ar@{=}[d]&\frac{\C[[p,q]]}{\ZN[[p,q]]+\tilde E^\Gamma_{\C,m+2}[[p]]+\tilde E^\Gamma_{\C,m+2}[[q]]+\C}\ar@{=}[d]\\\tilde W_{m+2}\ar[r]^{\tilde \pi}&W_{m+2}}$$
(see the end of the proof of \cite[Prop. 3.3.2]{laures-expansion}) shows that
the map $\tilde f=i\circ F$ factors as
\[ \tilde f:F^2\pi^S_m/F^3\pi^S_m\stackrel{f}{\to}\tilde{W}_{m+2}\stackrel{\tilde{\pi}}{\to}W_{m+2}\ ,\]
and $f$ is still injective after inverting $N$.

\begin{ddd}\label{deff}
We will call the  map $f:F^2\pi_m^S/F^3\pi_m^S\to \tilde{W}_{m+2}$ the $f$-invariant.
\end{ddd}
This map will be the basic object linking the analytical and 
topological indices $\eta^{an}$ and $\eta^{top}$
defined in \ref{etatop} and \ref{etaan}.

%
%
%
%

\section{The relation between  $\eta^{an}$ and $f$}\label{ufwiefwefwefwef}

In this Section we find the precise relation between $\eta^{an}$ and the $f$-invariant of Laures. The argument is based on Laures' geometric description of the $f$-invariant in terms of manifolds with corners (recalled in the preceeding section) and the Atiyah-Patodi-Singer type index theorem for manifolds with corners  \cite{math.DG/0201112}.
As a side result we get an analytic proof for the fact already known to Laures that
the $f$-invariant actually takes values in a  very small subgroup of $W_{m+2}$, cf. equation (\ref{dwsdasdaaw}). The properties of $\eta^{an}$ claimed in Theorem \ref{etaan}
can now be shown as a consequence of the known properties of the $f$-invariant. In Section \ref{idqwdqwdqwdwd} we will give independent analytic proofs for most of them.



We resume notations and assumptions as in Section \ref{uiewdwedwed}.
We choose  a Riemannian metric $g^{TX}$ on $X$ which is compatible with the corner structure. More precisely we assume that it is admissible in the sense of \cite{math.DG/0201112}, i.e. that we assume product structures near the boundary  components $Y_0,Y_1$ which meet with a right angle at the corner $Y_0\cap Y_1=Z$. The admissible Riemannian metric on $X$  gives rise to a Levi-Civita connection $\nabla^{LC}$. We further choose an extension $\nabla^{LC,L}$ of the Levi-Civita connection to a $Spin^c$-connection.

From now on we will distinguish the tangent bundle $TX$ from its stabilization $TX^s\cong TX\oplus (X\times \R^r)$.
We will further assume a metric on $TX^s$ 
such  that the decomposition $TX^s\cong T^0X\oplus T^1X$ is orthogonal,
the complex structures on $T^iX$ are anti-selfadjoint,
and such that that the induced metric on $T^iX_{|Y_{1-i}}$ is the metric given by the framing. Finally we assume a connection $\nabla^{T^iX}$ which preserves the splitting, the metric and the complex structure and restricts to the trivial connections on
$T^iX|_{Y_{1-i}}$. Note that the Levi-Civita connection can be extended by the trivial connection to a connection $\nabla^{LC,X}$ on $TX^s$ (which of course does not necessarily preserve the splitting or the complex structure).

%

 We abbreviate $$W(p,q):=\ch(\nabla^{C(T^0X)(p)})\wedge  \ch(\nabla^{C(T^1X)(q)})\in \Omega(X)\otimes E^\Gamma_\C[[p]]\otimes E^\Gamma_\C[[q]] \subset\Omega(X)[[p,q]]\ .$$
In the first step we replace $\Td(\nabla^{TX})$ by
$\Td(\nabla^{LC,L})$.
By Stoke's theorem we have
\begin{eqnarray}
\hat F(X)&=& \int_X \Td(\nabla^{T^0X} ) \wedge \Td(\nabla^{T^1X}) \wedge W(p,q)\nonumber\\
&=&\int_X \Td(\nabla^{LC,L}) W(p,q)+\int_X d \tilde \Td(\nabla^{T^0X}\oplus \nabla^{T^1X},\nabla^{LC,L})W(p,q)\nonumber\\
&=&
\int_X \Td(\nabla^{LC,L}) W(p,q)+\int_Y \tilde \Td(\nabla^{T^0X}\oplus \nabla^{T^1X},\nabla^{LC,L})W(p,q)\label{fzweufwef}\
,\end{eqnarray}
where $Y:=Y_0\cup Y_1$, and
$\tilde \Td(\nabla^{T^0X}\oplus \nabla^{T^1X},\nabla^{LC,L})$ is the transgression Todd form satisfying
$$d \tilde \Td(\nabla^{T^0X}\oplus \nabla^{T^1X},\nabla^{LC,L})=\Td(\nabla^{T^0X}\oplus \nabla^{T^1X})-\Td(\nabla^{LC,L})\ .$$
We can further write
\begin{eqnarray*}
\int_Y \tilde \Td(\nabla^{T^0X}\oplus \nabla^{T^1X},\nabla^{LC,L})W(p,q)&=&
\int_{Y_0} \tilde \Td(\nabla^{T^0X}\oplus \nabla^{T^1X},\nabla^{LC,L})W(p,q)\\&&+
\int_{Y_1} \tilde \Td(\nabla^{T^0X}\oplus \nabla^{T^1X},\nabla^{LC,L})W(p,q)\ .
\end{eqnarray*}
Since $\nabla^{T^iX}_{|Y_{1-i}}$ is trivial we have
$\ch(\nabla^{C(T^iX)(p)})_{|Y_{1-i}}=1$ and therefore
$$W(p,q)_{|Y_0}\in E^\Gamma_\C[[p]]\otimes \Omega(Y_0)\ ,\quad W(p,q)_{|Y_1}\in E^\Gamma_\C[[q]]\otimes \Omega(Y_1) .$$
Hence 
\begin{eqnarray}\int_{Y_0} \tilde \Td(\nabla^{T^0X}\oplus \nabla^{T^1X},\nabla^{LC,L})W(p,q)&\in& \C[[p]]\label{udiqwdqwd333}\\\int_{Y_1} \tilde \Td(\nabla^{T^0X}\oplus \nabla^{T^1X},\nabla^{LC,L})W(p,q)&\in& \C[[q]]\nonumber \ .
 \end{eqnarray}
Note that $\hat F(X)\in (E^\Gamma_\C\otimes E^\Gamma_\C)_{m+2}[[p,q]]$ while
the two terms on the right-hand side of (\ref{fzweufwef}) separately are inhomogeneous elements of 
$E^\Gamma_\C\otimes E^\Gamma_\C$.

We now can use the index theorem in order to express $F(X)$ in terms of the
$\partial <2>$-manifold $Y$. 
We assume that $m:=\dim(Z)>0$ is even.
We will ultimately look at the index of the twisted Dirac operator
$$\Dirac_X\otimes C(T^0X)(p)\otimes C(T^1X)(q)\ .$$ 
In order to turn this operator on a manifold with corners into a Fredholm operator we will choose a boundary taming.
Here we use the language introduced in \cite{math.DG/0201112}. The idea is to attach cylinders
to all boundary components and to  complete the corner by a quadrant so that we get a complete manifold with a Dirac type operator which is translation invariant at infinity. In order to turn this operator into a Fredholm operator we add smoothing perturbations to the operators on the boundary and corner faces to make them invertible. The notion of a boundary taming subsumes these choices.

In general there are obstructions to choosing a boundary taming but in the present case boundary tamings exist:

First of all, the operator
$\Dirac_Z$ bounds (actually in two ways through $Y_i$, $i=0,1$), and therefore
$\ind (\Dirac_Z)=0$. Hence it admits a taming
$\Dirac_{Z,t}$. Since
$$ [C(T^0X)(p)\otimes C(T^1X)(q)]_{|Z}$$
is a power series of trivial bundles we get an induced taming of
$$\Dirac_{Z,t}\otimes C(T^0X)(p)\otimes C(T^1X)(q)\ .$$
We interpret this choice as boundary tamings 
 $$(\Dirac_{Y_i}\otimes C(T^0X)(p)\otimes C(T^1X)(q))_{bt}$$
 of the faces $Y_i$.
 We can now extend these boundary tamings to tamings of
the faces $$(\Dirac_{Y_i}\otimes C(T^0X)(p)\otimes C(T^1X)(q))_t$$ since the manifolds $Y_i$ are odd-dimensional.
These choices make up the boundary taming
$$(\Dirac_X\otimes C(T^0X)(p)\otimes C(T^1X)(q))_{bt}\ .$$

The index theorem for manifolds with corners \cite{math.DG/0201112}  now gives
\begin{eqnarray}\lefteqn{
\int_X \Td(\nabla^{LC,L}) W(p,q)}&&\nonumber\\&&+\eta((\Dirac_{Y_0}\otimes C(T^0X)(p)\otimes C(T^1X)(q))_t)\nonumber\\&&+\eta((\Dirac_{Y_1}\otimes C(T^0X)(p)\otimes C(T^1X)(q))_t)\nonumber\\&&\hspace{2cm}=\ind((\Dirac_X\otimes C(T^0X)(p)\otimes C(T^1X)(q))_{bt})\label{dwqdqwdwqd}\\&\in&\ZN[[p,q]]\nonumber\ .
\end{eqnarray}
If we combine (\ref{fzweufwef}) and (\ref{dwqdqwdwqd}), then we get
an equality in $W_{m+2}$
\begin{eqnarray}\label{1234}
\lefteqn{f(X)}&&\\&\hspace{-0.5cm}=&\hspace{-0.5cm}\int_{Y_0} \tilde \Td(\nabla^{T^0X}\oplus \nabla^{T^1X},\nabla^{LC,L})W(p,q)-\eta((\Dirac_{Y_0}\otimes C(T^0X)(p)\otimes C(T^1X)(q))_t)\label{teqzweqwe1}\\&\hspace{-0.5cm}+&\hspace{-0.5cm}\int_{Y_1} \tilde \Td(\nabla^{T^0X}\oplus \nabla^{T^1X},\nabla^{LC,L})W(p,q)-\eta((\Dirac_{Y_1}\otimes C(T^0X)(p)\otimes C(T^1X)(q))_t)\label{teqzweqwe2}\ .\end{eqnarray}
Let us now consider the first term associated to $Y_0$. Since $T^1Y_0$ is trivial we see that $(D_{Y_0}\otimes C(T^0X)(p)\otimes C(T^1X)(q))_{bt}$
is a sum of copies of $(\Dirac_{Y_0}\otimes C(T^0X)(p))_{bt}$.
We first choose an extension of this boundary taming to a taming and then let $(\Dirac_{Y_0}\otimes C(T^0X)(p)\otimes C(T^1X)(q))_{t}$ be the induced taming. With these choices we have
$$\eta((\Dirac_{Y_0}\otimes C(T^0X)(p)\otimes C(T^1X)(q))_{t})\in \C[[p]]\ .$$ 
By a similar choice we ensure that
$$\eta((\Dirac_{Y_1}\otimes C(T^0X)(p)\otimes C(T^1X)(q))_{t})\in \C[[q]]\ .$$ 
 
Using (\ref{1234}) we conclude that
$$\hat F(X)\in (\ZN[[p,q]]+\C[[p]]+\C[[q]])\cap (E^\Gamma_\C\otimes E^\Gamma_\C)_{m+2}[[p,q]]\ .$$

Let us consider the subgroup
$$U_{m+2}:=\frac{\C[[p]]+\C[[q]]}{\ZN[[p]]+\ZN[[q]]+E^\Gamma_{\C,m+2}[[p]]+E_{\C,m+2}^\Gamma[[q]]+\C}\subseteq \tilde{W}_{m+2}\ .$$

We can split
$$U_{m+2}=U_{m+2}^p\oplus U_{m+2}^q$$ with
$$U_{m+2}^p:=\frac{\C[[p]]}{\ZN[[p]]+E^\Gamma_{\C,m+2}[[p]]+\C}\ ,\quad U_{m+2}^q:=\frac{\C[[q]]}{\ZN[[q]]+E^\Gamma_{\C,m+2}[[q]]+\C}.$$
These are exactly the groups where the analytical index $\eta^{an}(Z)$ lives.
We see that $f(Z)=i(F(X))$ is represented by a pair 
$$\tilde f(Y_0)\oplus \tilde f(Y_1)\in U_{m+2}^p\oplus U^q_{m+2}\ ,$$ where
$$\tilde f(Y_0):=[(\ref{teqzweqwe1})]\ ,\quad \tilde f(Y_1):=[(\ref{teqzweqwe2})]\ ,$$
and the brackets $[\dots]$ mean that we take the classes of the formal power series in the corresponding quotient
$U^q_{m+2}$ or $U^p_{m+2}$, respectively.

Using the fact that
$T^1Y_0$ is trivialized we can simplify the expression for $\tilde f(Y_0)$ further. 
We get  
\begin{eqnarray*}
\tilde f(Y_0)&=&[\int_{Y_0} \tilde \Td(\nabla^{T^0Y_0}\oplus \nabla^{T^1Y_0},\nabla^{LC,L})\wedge \ch(\nabla^{C(T^0Y_0)(p)})
-\eta((\Dirac_{Y_0}\otimes C(T^0Y_0)(p))_t)]\\
&=&[\int_{Y_0} \tilde \Td(\nabla^{TY_0},\nabla^{LC,L})\wedge \ch(\nabla^{C(TY_0)(p)})
-\eta((\Dirac_{Y_0}\otimes C(TY_0)(p))_t)]
\\&=&\eta^{an}(Z)\ .
\end{eqnarray*}

In a similar way we get
$$\tilde f(Y_1)=-\eta^{an}(Z)\ ,$$
where the sign arises since we orient $Z$ as the boundary of $Y_0$, and this orientation is opposite to the orientation of $Z$ as the boundary of $Y_1$.


Combining the above, we obtain
\begin{equation}\label{dwsdasdaaw}f(Z)=\eta^{an}(Z)(p)\oplus -\eta^{an}(Z)(q).\end{equation}

The prescription $q\mapsto 0$ induces a projection \begin{equation}\label{duwqidwqd}\pi:\frac{\C [[p,q]]}{\ZN[[p,q]]+ E^\Gamma_{\C,m+2}[[q]]+ E^\Gamma_{\C,m+2}[[p]]+\C}\to \frac{\C [[p]]}{\ZN[[p]]+ E^\Gamma_{\C,m+2}[[p]]+\C}\ ,\end{equation}  i.e. a map
$\pi:\tilde{W}_{m+2}\to U^p_{m+2}$. We get
$\eta^{an}(Z)(p)=\pi(f(Z))$ in $U^p_{m+2}$.\\

\begin{proof} (of Theorem \ref{etaan})\\
From the above we have a commutative diagram

\[ \xymatrix{ F^2\pi^S_m/F^3\pi_m^S\otimes_\Z \Z\left[\frac{1}{N}\right]\ar[r]^(.7)f & \tilde{W}_{m+2} \ar[r]^\pi & U_{m+2}^q\\ F^2\pi_m^S\ar[u]\ar[rru]_{\eta^{an}}} \]

and the composition $\pi\circ f$ is injective according to \cite[Lemma 3.2.2]{laures-expansion}.
\end{proof}

\section{The relation between $\eta^{top}$ and $f$}\label{csacoascascs}

Let $m\ge 2$ be even and $\alpha\in \pi_m^S$.
Recall that $f(\alpha)\in \tilde{W}_{m+2}$ and $\eta^{top}(\alpha)\in U^p_{m+2}$, and that we have introduced a map $\pi:\tilde{W}_{m+2}\to U^p_{m+2}$ above, see (\ref{duwqidwqd}).
\begin{prop}\label{uigfvwegffwef}
We have $\pi(f(\alpha))=\eta^{top}(\alpha)$.
\end{prop}
\proof
We resume notation and assumptions from the Adams resolution (\ref{ufwiefewfwef}) and
consider the following web of horizontal and vertical fiber sequences
constructed by suitably smashing the defining fiber sequences
$$ \xymatrix{\overline{MU}\ar[r]&\overline{MU}_\Q\ar[r]&\overline{MU}_{\Q/Z}}$$
and
$$\xymatrix{ S\ar[r]& MU\ar[r] & \overline{MU}\ .}$$

%
\begin{equation}\label{web1}
\xymatrix{
\Sigma^{-1}\overline{MU}_\Q\wedge MU\ar[r]\ar[d]&\Sigma^{-1}\overline{MU}_{\Q/\Z}\wedge MU \ar[r]\ar[d]&\overline{MU}\wedge MU\ar[r]\ar[d]& \overline{MU}_\Q\wedge MU\ar[d]\\\Sigma^{-1}\overline{MU}_\Q\wedge \overline{MU}\ar[r]\ar[d]&\Sigma^{-1}\overline{MU}_{\Q/\Z}\wedge \overline{MU}\ar[d]\ar[r]&\overline{MU}\wedge \overline{MU}\ar[d] \ar[r]&  \overline{MU}_\Q\wedge \overline{MU} \ar[d] \\\overline{MU}_\Q\ar[d]\ar[r]&\overline{MU}_{\Q/\Z}\ar[d]\ar[r]&\Sigma\overline{MU}\ar[d]\ar[r]&\Sigma \overline{MU}_\Q\ar[d]\\
\overline{MU}_\Q\wedge MU\ar[r]\ar[d]&\overline{MU}_{\Q/\Z}\wedge MU \ar[r]\ar[d]&\Sigma\overline{MU}\wedge MU\ar[r]\ar[d]&\Sigma \overline{MU}_\Q\wedge MU\ar[d]\\\overline{MU}_\Q\wedge \overline{MU}\ar[r] &\overline{MU}_{\Q/\Z}\wedge \overline{MU} \ar[r]&\Sigma\overline{MU}\wedge \overline{MU}  \ar[r]&\Sigma \overline{MU}_\Q\wedge \overline{MU}.  }
\end{equation}

The class $\hat \alpha\in \overline{MU}_{m+1}$ is torsion and therefore has a lift $\tilde \alpha_{\Q/\Z}\in \overline{MU}_{\Q/\Z,m+2}$. Since $\hat \alpha$ 
admits the lift $\tilde{\alpha}$ in (\ref{ufwiefewfwef}), it is in the 
kernel of $\id\wedge\epsilon: \overline{MU}_{m+1}\to (\overline{MU}\wedge MU)_{m+1}$, hence the image of $\tilde \alpha_{\Q/\Z}$ under
$\overline{MU}_{\Q/\Z,m+2}\to (\overline{MU}_{\Q/\Z}\wedge MU)_{m+2}$ further lifts to some $\overline{\eta}\in (\overline{MU}_\Q\wedge MU)_{m+2}$, c.f. (\ref{diag2}).
The image of $\overline{\eta}$ under the map
$$\bar\nu_\Q\wedge \theta: \overline{MU}_\Q\wedge MU \to \overline{T}_\Q\wedge K\, ,\, \overline{\nu}_\Q:=(\overline{\gamma}\circ\overline{\kappa})_\Q$$
is a possible choice of the  element  $\eta\in (\bar T_\Q\wedge K)_{m+2}$ in the construction of $\eta^{top}$, c.f. (\ref{diag2}).
By a diagram chase one checks that the class
$\overline{\eta}$ projects under
$$\overline{MU}_\Q\wedge MU\to \overline{MU}_\Q\wedge \overline{MU}$$
 to the image $-\tilde \alpha_\Q$ of the element $-\tilde \alpha\in (\overline{MU}\wedge \overline{MU})_{m+2}$ from (\ref{ufwiefewfwef}) under the map
$$\overline{MU}\wedge \overline{MU}\to \overline{MU}_\Q \wedge\overline{MU}\ .$$

We summarize the above discussion in the following diagram.

$$\xymatrix{\overline{\eta}\ar@{|.>}[dd]^{(\ref{diag2})}\ar@{|.>}[rrr]^{\mbox{(diagram chase)}} &&&-\tilde \alpha_\Q\ar@{|.>}[dddd]^{(Def. \ref{deff})}\\
& (\overline{MU}_\Q \wedge MU)_{m+2}\ar[r]\ar[d]^{\bar \nu_\Q \wedge \theta    }&(\overline{MU}_\Q\wedge \overline{MU})_{m+2}\ar[dd]^{(p,q)-\mbox{\small expansion}\circ \bar \kappa\wedge \bar\kappa}&\\
\eta\ar@{|.>}[dd]^{(Def. \ref{etatop})}& 
(\bar T_\Q \wedge K)_{m+2}\ar[d]&&\\
&\frac{\C [[p]]}{\ZN[[p]]+ E^\Gamma_{\C,m+2}[[p]]+\C}& \frac{\C [[p,q]]}{\ZN[[p,q]]+ E^\Gamma_{\C,m+2}[[q]]+ E^\Gamma_{\C,m+2}[[p]]+\C}\ar[l]_(.6){\pi }&\\
-\eta^{top}(\alpha)=-\pi(f(\alpha))&&&-f(\alpha).\ar@{|.>}[lll]}
$$

Mapping $\overline{\eta}$ clockwise to $U_{m+2}^p$ yields $-\pi(f(\alpha))$
while mapping it counter-clockwise gives $-\eta^{top}(\alpha)$. We claim that
the solid diagram above commutes. This immediately implies that
$\eta^{top}(\alpha)=\pi(f(\alpha))$.
In order to see the claim note that we can factorize the orientation $\theta:MU\to K$ as
$$MU\stackrel{\kappa}{\to} \tilde{E}^\Gamma\stackrel{\gamma}{\to} T\stackrel{q\mapsto 0}{\to} K\ .$$
This is applied to the second factor.\hB

\section{Analysis of $\eta^{an}$}\label{idqwdqwdqwdwd}

\label{xsjxkasxasx}\label{dusadasdasd}


In this section we present
the construction of $\eta^{an}$ in complete generality which requires the use of tamings.
In Theorem \ref{main} we show using the arguments of an index theorist that
$\eta^{an}$ is independent of a plethora of auxiliary choices
and factors over the framed bordism  group, thus reproving parts of Theorem \ref{idqwdqwdqwdwd}. Finally we prove the tertiary index Theorem \ref{dsadioasdlsadasdsadsadd}.

Let $m>0$ be even and
assume that the class $\alpha\in \pi^S_m\cong \Omega^{fr}_m$ is represented by a manifold $Z$ with a framing of the stable tangent bundle $TZ^s$.   Since $\alpha\in \pi^S_m$ is a torsion element, and $MU_m$ is torsion-free, the image $\epsilon(\alpha)\in MU_m$ under the unit $\epsilon:S\to MU$ vanishes.
Hence we can choose a zero bordism $N$, $\partial N\cong Z$, with a stable complex structure on $TN^s$ which  extends the framing.



 We choose a Riemannian metric on $N$ with a product structure which induces a Riemannian metric on $Z$.  We choose furthermore a hermitian metric and a hermitian connection on $TN^s$ which become the trivial ones near $Z$.



The normal complex structures on $N$ and $Z$ determine a $Spin^c$-structure.
We choose an extension of the Levi-Civita connection $\nabla^{LC}$ on $N$ to a $Spin^c$-connection (see Section \ref{duwiqdqwd})
which is of product type near $Z$. With the complex spinor bundle,
$N$ becomes a geometric manifold $\cN $ with boundary $\cZ=\partial \cN $.
We refer to \cite{math.DG/0201112} for the notion of a geometric manifold which is used as a shorthand for the collection of structures needed to define a generalized Dirac operator $\Dirac_N$. The relation $\cZ=\partial \cN$ 
implies that the boundary reduction of
$\Dirac_N$ is $\Dirac_Z$.

It follows from the bordism invariance of the index that
$\ind(\Dirac_Z )=0$. Therefore we can choose some taming $\Dirac_{Z,t}$ (see Section \ref{ufwiefwefwefwef} and \cite{math.DG/0201112}). For the present paper it suffices to understand that a taming is a choice of smoothing operators on all faces of a manifold with corners $M$ which can be used to turn the Dirac operator $\Dirac_{M}$ into a Fredholm operator  $\Dirac_{M,t}$ to which the methods of local index theory apply.
Note that in the present note we
use a different notation which attaches the taming to the symbol for Dirac operator instead of the geometric manifold. The operator $\Dirac_{Z,t}$ is thus an invertible
perturbation of $\Dirac_Z$. If the latter itself is invertible, then the trivial taming is a canonical choice
used in Section \ref{sdvdsv}.

Recall the definition (\ref{euwdioewdwed}) of the bundles $W_n\to N$ as coefficients of the formal power series $C(TN^s)(p)$. These bundles come with induced hermitian metrics and hermitian connections $\nabla^{W_n}$.
The trivialization of $TN^s$ near $Z$ induces trivializations of $W_n$ near $Z$.
Hence we have identifications of $\Dirac_Z\otimes W_{n|Z}$ with direct sums of copies of $\Dirac_Z$. We see that the taming $\Dirac_{Z,t}$ induces a boundary taming
$(\Dirac_N\otimes W_n)_{bt}$.

%

Since $N$ is odd-dimensional we can extend this boundary taming to a taming $(D_N\otimes W_n)_{t}$.
The sequence of $\eta$-invariants
$\eta((\Dirac_N \otimes W_n)_t)\in \R$ gives rise to a formal power series which we will denote by (compare (\ref{fuwehiewfewf}))
\begin{equation}\label{dioqwdqwd}\eta(p):= \eta((\Dirac_N\otimes C(TN^s)(p))_t)\in \C[[p]]\ .\end{equation}

\begin{ddd}\label{etwer2a}
We define
$$\eta^{an}\in  \frac{\C[[p]]}{\ZN[[p]]+E^\Gamma_{\C,m+2}[[p]]+\C}$$
as the class represented by 
$$\int_N \tilde \Td(\nabla^{TN^s},\nabla^{LC,L})\wedge \ch(\nabla^{C(TN^s)(p)})-\eta(p)\ .$$
\end{ddd}

\begin{theorem}\label{main}\label{qwduiqwdqwdqwd}
The element $\eta^{an}$ does only depend on the class $\alpha\in \pi^S_m$.
\end{theorem}
Since $\eta^{an}$ is clearly additive under disjoint union of framed manifolds and changes sign if we  
switch the orientation we thus  get a homomorphism
$$\eta^{an}:\pi^S_m\to   \frac{\C[[p]]}{\ZN[[p]]+E^\Gamma_{\C,m+2}[[p]]+\C}\ .$$
We first show the independence of $\eta^{an}$
of the various choices in the construction.

\begin{lem}
The class $\eta^{an}$ does not depend on the choice of the extension
$(\Dirac\otimes C(TN^s)(p))_t$ of the boundary taming.
\end{lem}
\proof
If $(\Dirac\otimes C(TN^s)(p))_t^\prime$ is a second choice with resulting $\eta^{\prime}(p)$ and $\eta^{an \prime}$, then by \cite[2.2.17]{math.DG/0201112}
$$\eta^\prime(p)-\eta(p)=\Sf((\Dirac\otimes C(TN^s)(p))_t^\prime,(\Dirac\otimes C(TN^s)(p))_t)\in \ZN[[p]]\ ,$$
where $\Sf(D_t,D_t^\prime)$ denotes the spectral flow of a family of pre-tamed Dirac operators interpolating between $D_t$ and $D_t^\prime$. This implies that $\eta^{an}=\eta^{an  \prime}$.
\hB

\begin{lem}
The class $\eta^{an}$ does not depend on the choice of the taming $\Dirac_{Z,t}$.
\end{lem}
\proof
Let $\Dirac_{Z,t}^\prime$ be a second choice.
We consider the product $\cZ\times I$. The two tamings $\Dirac_{Z,t}, \Dirac_{Z,t}^\prime$ induce a boundary taming $\Dirac_{Z\times I,bt}$. This boundary taming 
can be extended to a taming
$\Dirac_{Z\times I,t}$  since $\cZ\times I$ is odd-dimensional.
The boundary of $N\times I$ consists of the faces $N\times \{0\}$, $N\times \{1\}$, and $Z\times I$.
We choose some extensions $(\Dirac_{N}\otimes C(TN^s))_t$, $(\Dirac_{N}\otimes C(TN^s))^\prime_t$   of the boundary tamings $\Dirac_{Z,t}\otimes C(TN^s)_{|Z}$ and $\Dirac_{Z,t}^\prime\otimes C(TN^s)_{|Z}$.
These choices give tamings of the
 the corresponding boundary face reductions of $(\Dirac_{N\times I}\otimes C(\pr_1^*TN^s))$. Together with the taming
$\Dirac_{Z\times I,t}\otimes C(TN^s)_{|Z}$ this yields a boundary taming
$(\Dirac_{N\times I}\otimes C(\pr_1^* TN^s))_{bt}$.
We now apply the index theorem \cite[Theorem 2.2.13 (2)]{math.DG/0201112} and get
\begin{eqnarray*}\lefteqn{\ind((\Dirac_{N\times I}\otimes C(\pr_1^* TN^s))_{bt})}&&\\&& =\eta(D_{\partial (N\times I)} \otimes C(TN^s)_{|\partial (N\times I)})_{bt} )+\Omega((\cN\times I)\otimes C(\pr_1^*TN^s))\in \ZN[[p]]\ ,\end{eqnarray*}
where $\eta(D_{\partial (N\times I)} \otimes C(\pr_1^*TN^s){|\partial (N\times I)})_{t}$
is the sum of the $\eta$-invariants of the boundary faces, i.e.
\begin{eqnarray*}
\eta(D_{\partial (N\times I)} \otimes C( \pr_1^*TN^s)_{|\partial (N\times I)})_{t}&=&
\eta(\Dirac_{Z\times I,t}^\prime\otimes C(\pr_1^*TN^s_{|Z}))\\&&-\eta((\Dirac_{N}\otimes C(TN^s))_t)\\&&+ \eta((\Dirac_{N}\otimes C(TN^s))^\prime_t)\ ,
\end{eqnarray*}
and $\Omega((\cN\times I)\otimes C(\pr_1^*TN^s))$ denotes the local contribution to the index.
Since the geometry of $(\cN\times I)$  is of product type we get
$\Omega((\cN\times I)\otimes C(\pr_1^* TN^s))=0$.
Furthermore, we have by (\ref{duqwdwqdqwd}) 
$$\eta(\Dirac_{Z\times I,t}^\prime\otimes C(\pr_1^*TN^s_{|Z})(p))\in \C\subset \C[[p]]\ ,$$
since $\pr_1^*TN^s_{|Z}$ is trivial.
This implies that $$\eta(\Dirac_{N}\otimes C(TN^s)(p))_t)\equiv\eta(\Dirac_{N}\otimes C(TN^s)(p))^\prime_t)\quad \mbox{modulo}\quad 
\ZN[[p]]+ \C$$  and hence the assertion of the Lemma. \hB

%
%

\begin{lem}
The class $\eta^{an}$ does not depend on the choice of the zero bordism $N$.
\end{lem}
\proof
Let $N^{\prime}$ be a second choice leading to $\eta^{an\prime}$.
Then we can form the closed manifold $Y:=N\cup_{Z}(N^\prime)^{op}$ by glueing $N$ and $N^\prime$ along their  boundaries.
We can choose the geometric structures on $N$ and $N^\prime$ (Riemannian metrics, $Spin^c$-connections and connections on stable tangent bundles) such that they coincide
near $Z$ and thus induce corresponding  geometric structures on $Y$. We let $\cY$ denote the corresponding geometric manifold.
Since $Y$ is odd-dimensional we can choose a taming $(\Dirac_Y\otimes C(TY^s))_t$.
The glueing formula for $\eta$-invariants gives
$$\eta((\Dirac_N\otimes C(TN^s)(p))_t)- \eta((\Dirac_{N^\prime}\otimes C(TN^{\prime s})(p))_t)- \eta((\Dirac_Y\otimes C(TY^s)(p))_t)\in \ZN[[p]]\ .$$
The calculation 
(\ref{udwqdqwdqwd}) together with the identity
\begin{eqnarray*}0&=&\int_N \tilde \Td(\nabla^{TN^s},\nabla^{LC,L})\wedge \ch(\nabla^{C(TN^s)(p)})-\int_{N^\prime} \tilde \Td(\nabla^{TN^{\prime s}},\nabla^{LC,L})\wedge \ch(\nabla^{C(TN^{\prime s})(p)})\\&&-\int_Y \tilde \Td(\nabla^{TY^s},\nabla^{LC,L})\wedge \ch(\nabla^{C(TY^s)(p)})=0\end{eqnarray*}
now implies that  
$\eta^{an}=\eta^{an\prime}$.
\hB



\begin{lem}
The class $\eta^{an}$ does only depend on the framed bordism class $\alpha$.
\end{lem}
\proof
Note that $\eta^{an}$ is additive with respect to disjoint union and changes sign if we reverse the orientation.
If $Z$ is framed zero bordant, then we can use this zero bordism in place of $N$. In this case
the bundle $TN^s$ is trivialized. We first extend the taming $\Dirac_{Z,t}$ to a taming $\Dirac_{N,t}$. It induces a taming
$\Dirac_{N,t}\otimes C(TN^s)$, and we get $$\int_N \tilde \Td(\nabla^{TN^s},\nabla^{LC,L})\wedge \ch(\nabla^{C(TN^s)(p)})-\eta((\Dirac_{N,t}\otimes C(TN^s)(p))_t)\in \C\ .$$
This implies the result. \hB
This finishes the proof of Theorem \ref{main}. \hB

 Recall the definition of
 $\eta^{top}$ given in Section \ref{wqdqwodqwdqwdqwd}.


\begin{theorem}\label{dwqdqwdwqdwd}
For even $m>0$ we have the equality of homomorphisms
$$\eta^{an}=\eta^{top}:\pi^S_m\to \frac{\C[[p]]}{\ZN[[p]]+E^\Gamma_{\C,m+2}[[p]]+\C}$$
\end{theorem}
\proof 
We apply Proposition \ref{uigfvwegffwef} to the equation (\ref{dwsdasdaaw}).   \hB

\section{Mod $k$-indices}\label{ddqwidwqudqwidwqd}

In the present Section we explain a way to represent
 $\eta^{an}$ as an analytic mod-$k$-index in the sense of Freed-Melrose \cite{MR1144425} or Higson \cite{MR1060635}.  This may open a different path to topological calculations of $\eta^{an}$, but note that at the moment the necessary generalization of the  $\Z/k\Z$-index theorem to manifolds with corners
is not available. 

Assume that $m>0$ is even and let
$\alpha\in \pi_m^S$ be represented by the stably framed manifold $Z$.  
Then there is a pair $(N,Z)$ consisting of the stably framed manifold $Z$ and a stably complex zero bordism $N$ which represents the class
$\hat \alpha\in \overline{MU}_{m+1}$ in (\ref{duqwidwqdwqdqd}).
We have seen that $\hat \alpha$ is a torsion class. Let $k>0$ be an integer such that  $k\hat \alpha=0$.
This means that there exists a manifold  $Y$ with corners  of codimension two and two boundary faces $\partial_iY$, $i=0,1$, and complex stable tangent bundle $TY^s\to Y$ such that
\begin{enumerate}
\item $\partial_0Y\cong k N$ as stably complex manifolds, where $k N$ is the disjoint union of $k$ copies of $N$,
\item the complex structure of $TY^s_{|\partial_1Y}$ refines to a framing,
\item the framing of $TY^s_{|k Z}$ is the given one on the $k$ copies of $Z$.  
\end{enumerate}
We choose the geometric structures (Riemannian metrics, $Spin^c$-connections and hermitian connections on the stable tangent bundles) adapted to the corner structure (as in Section \ref{ufwiefwefwefwef}) and get a geometric manifold
$\cY$ so that $\partial_0\cY=k\cN$. We extend the taming $\Dirac_{kZ,t}$ (which is induced by $\Dirac_{Z,t}$) to a taming  $\Dirac_{\partial_1 Y,t}$
(this is possible since this boundary is odd-dimensional).
It induces a taming $\Dirac_{\partial_1 Y,t} \otimes C(TY^s_{|\partial_1 Y})$.
Together with a taming $(\Dirac_{\partial_0Y}\otimes C(TY^s_{|\partial_0Y}))_t$ induced by $k$ copies of the  taming $(\Dirac_N\otimes C(TN^s))_t$
 this yields a boundary taming
$(\Dirac_Y\otimes C(TY^s))_{bt}$.
\begin{prop} In the above situation we have
$$\eta^{an}(\alpha)=[-\frac{1}{k}\ind((\Dirac_Y\otimes C(TY^s)(p))_{bt})] \in \frac{\C[[p]]}{\ZN[[p]]+E^\Gamma_{\C,m+2}[[p]]+\C}\ .$$
\end{prop}
\proof
We have the index theorem for manifolds with corners \cite[Theorem 2.2.13 (2)]{math.DG/0201112}
\begin{eqnarray}
\lefteqn{\ind((\Dirac_Y\otimes C(TY^s)(p))_{bt})}&&\label{dqwidoqwdqwd}\\&=&\Omega(\cY\otimes C(TY^s)(p))+ \eta(\Dirac_{\partial_1 Y,t}\otimes C(TY^s_{|\partial_1Y})(p))\nonumber\\&&+k\eta((\Dirac_N\otimes C(TY^s_{|\partial_0 Y})(p))_t)\nonumber\\&\in&\ZN[[p]]\ . 
\end{eqnarray}
We now observe that  $\eta(\Dirac_{\partial_1 Y,t}\otimes C(TY^s_{|\partial_1Y})(p))\in \C$,  
and
\begin{eqnarray*}
\Omega(\cY\otimes C(TY^s)(p))&=&
\int_Y \Td(\nabla^{LC,L})\wedge \ch(\nabla^{C(TY^s)(p)})\\
&=&\int_Y \Td(\nabla^{TY^s})\wedge \ch(\nabla^{C(TY^s)(p)})\\&&-\int_{\partial Y} \tilde\Td(\nabla^{TY^s}, \nabla^{LC,L})\wedge \ch(\nabla^{C(TY^s)(p)})
\end{eqnarray*}
(compare (\ref{uidqwdqwd})). The latter equality
 shows that
$$\Omega(\cY\otimes C(TY^s)(p))+\int_{\partial Y} \tilde\Td(\nabla^{TY^s}, \nabla^{LC,L})\wedge \ch(\nabla^{C(TY^s)(p)})\in E^\Gamma_{\C,m+2}[[p]]\ .$$
We further observe that
$$\int_{\partial_1Y}\tilde\Td(\nabla^{TY^s}, \nabla^{LC,L})\wedge \ch(\nabla^{C(TY^s)(p)})\in \C$$
and 
$$\int_{\partial_0Y}\tilde\Td(\nabla^{TY^s}, \nabla^{LC,L})\wedge \ch(\nabla^{C(TY^s)(p)})=k \int_{N}\tilde\Td(\nabla^{TN^s}, \nabla^{LC,L})\wedge \ch(\nabla^{C(TN^s)(p)})\ .$$

We conclude that
\begin{eqnarray*}\lefteqn{\ind((\Dirac_Y\otimes C(TY^s)(p))_{bt})}&&\\&\equiv& k\eta((\Dirac_N\otimes C(TY^s_{|\partial_0 Y})(p))_t)\\&&-k\int_{N}\tilde\Td(\nabla^{LC,L}, \nabla^{TN^s})\wedge \ch(\nabla^{C(TN^s)(p)})\\&=&-k \eta(p)\end{eqnarray*}
modulo $E^\Gamma_{\C,m+2}[[p]]+\C$,
where $\eta(p)$ is as in (\ref{dioqwdqwd}).
Therefore
$$\eta^{an}(\alpha)=[-\frac{1}{k}\ind((\Dirac_Y\otimes C(TY^s)(p))_{bt})] \in \frac{\C[[p]]}{\ZN[[p]]+E^\Gamma_{\C,m+2}[[p]]+\C}\ .$$ \hB


\begin{thebibliography}{APS75b}

\bibitem[Ada66]{adams}
J.~F. Adams.
\newblock On the groups {$J(X)$}. {IV}.
\newblock {\em Topology}, 5:21--71, 1966.

\bibitem[AHS01]{ahs}
M.~Ando, M.~J. Hopkins, and N.~P. Strickland.
\newblock Elliptic spectra, the {W}itten genus and the theorem of the cube.
\newblock {\em Invent. Math.}, 146(3):595--687, 2001.

\bibitem[AS68]{MR0236950}
M.~F. Atiyah and I.~M. Singer.
\newblock The index of elliptic operators. {I}.
\newblock {\em Ann. of Math. (2)}, 87:484--530, 1968.

\bibitem[APS75a]{MR0397797}
M.~F. Atiyah, V.~K. Patodi, and I.~M. Singer.
\newblock Spectral asymmetry and {R}iemannian geometry. {I}.
\newblock {\em Math. Proc. Cambridge Philos. Soc.}, 77:43--69, 1975.

\bibitem[APS75b]{MR0397798}
M.~F. Atiyah, V.~K. Patodi, and I.~M. Singer.
\newblock Spectral asymmetry and {R}iemannian geometry. {II}.
\newblock {\em Math. Proc. Cambridge Philos. Soc.}, 78(3):405--432, 1975.

\bibitem[AS69]{MR0259946}
M.~F. Atiyah and G.~B. Segal.
\newblock Equivariant {$K$}-theory and completion.
\newblock {\em J. Differential Geometry}, 3:1--18, 1969.

\bibitem[BL]{mark2line}
M. Behrens, G. Laures.
\newblock Beta-family congruences and the f-invariant.
\newblock {\em to appear: Geometry and Topology Monographs}, available at:
http://front.math.ucdavis.edu/0809.1125


\bibitem[BP04]{MR2085661}
U. Bunke and J. Park.
\newblock Determinant bundles, boundaries, and surgery.
\newblock {\em J. Geom. Phys.}, 52(1):28--43, 2004.

\bibitem[B09]{math.DG/0201112}
U.~Bunke.
\newblock Index theory, eta forms, and {D}eligne cohomology.
\newblock {\em Mem. Amer. Math. Soc.}, 198(928):vi+120, 2009.

\bibitem[BS]{bunke-2007}
U. Bunke, T. Schick.
\newblock Smooth $K$-theory.
\newblock Preprint 2007, accepted. ArXiv.org:0905.4181

\bibitem[Bod]{bod} 
H. von Bodecker
\newblock On the geometry of the f-invariant. 
\newblock Preprint 2008, http://front.math.ucdavis.edu/0808.0428

\bibitem[DS84]{deninger}
C. Deninger and W. Singhof.
\newblock The $e$-invariant and the spectrum of the Laplacian for compact nilmanifolds covered by Heisenberg groups. 
\newblock {\em  Invent. Math.}, 78  (1984),  no. 1, 101--112. 

\bibitem[FM92]{MR1144425}
D. S. Freed and R. B. Melrose.
\newblock A mod {$k$} index theorem.
\newblock {\em Invent. Math.}, 107(2):283--299, 1992.

\bibitem[Fra92]{MR1235295}
J. Franke.
\newblock On the construction of elliptic cohomology.
\newblock {\em Math. Nachr.}, 158:43--65, 1992.

\bibitem[GHMR05]{k2local}
P. Goerss, H.-W. Henn, M. Mahowald and C. Rezk.
\newblock A resolution of the $K(2)$-local sphere at the prime 3.
\newblock {\em  Ann. of Math.} (2), 162  (2005),  no. 2, 777--822. 


\bibitem[Hig90]{MR1060635}
N.~Higson.
\newblock An approach to {${\bf Z}/k$}-index theory.
\newblock {\em Internat. J. Math.}, 1(2):189--210, 1990.

\bibitem[HHR]{hhr}
A. Hill, M. J. Hopkins, D. C. Ravenel
\newblock On the non-existence of elements of Kervaire invariant one
\newblock  Preprint 2009, http://arxiv.org/abs/0908.3724

\bibitem[HBJ92]{hirzebruch}
F. Hirzebruch, T. Berger, and R. Jung.
\newblock {\em Manifolds and modular forms}.
\newblock Aspects of Mathematics, E20. Friedr. Vieweg \& Sohn, Braunschweig,
  1992.
\newblock With appendices by Nils-Peter Skoruppa and by Paul Baum.

\bibitem[HN07]{hornbostelnaumann}
J. Hornbostel and N. Naumann.
\newblock Beta-elements and divided congruences.
\newblock {\em Amer. J. Math.}, 129(5):1377--1402, 2007.

\bibitem[Lau99]{laures-expansion}
G. Laures.
\newblock The topological {$q$}-expansion principle.
\newblock {\em Topology}, 38(2):387--425, 1999.

\bibitem[Lau00]{laures-cob}
G. Laures.
\newblock On cobordism of manifolds with corners.
\newblock {\em Trans. Amer. Math. Soc.}, 352(12):5667--5688 (electronic), 2000.

\bibitem[MRW77]{mrw}
H. Miller, D. Ravenel, S. Wilson.
\newblock Periodic phenomena in the Adams-Novikov spectral sequence. 
\newblock {\em Ann. Math.} (2), 106 (1977), no. 3, 469--516.


\bibitem[Rav86]{ravenel}
D. Ravenel.
\newblock {\em Complex cobordism and stable homotopy groups of spheres}, volume
  121 of {\em Pure and Applied Mathematics}.
\newblock Academic Press Inc., Orlando, FL, 1986.


\bibitem[Shi71]{shimura}
G. Shimura.
\newblock {\em Introduction to the arithmetic theory of automorphic
              functions}, Kan{\^o} Memorial Lectures, No. 1
\newblock Publications of the Mathematical Society of Japan, No. 11.
              Iwanami Shoten, Publishers, Tokyo

\end{thebibliography}
\end{document}